\documentclass[smallextended,envcountsect]{svjour3}
\smartqed 
\usepackage{graphicx}
\journalname{JOTA}

\usepackage{amssymb}
\usepackage{amsmath}

\usepackage{caption}
\usepackage{cite}
\usepackage{float}
\floatstyle{plaintop}
\restylefloat{table}

\usepackage{hyperref}

\def\Def{:=}

\def\beq{\begin{equation}}
\def\eeq{\end{equation}}

\def\R{\mathbb{R}}
\def\E{\mathbb{E}}

\def\BI{\begin{itemize}}
\def\EI{\end{itemize}}

\newcommand{\refLE}[1]{\ensuremath{\stackrel{(\ref{#1})}{\leq}}}
\newcommand{\refEQ}[1]{\ensuremath{\stackrel{(\ref{#1})}{=}}}
\newcommand{\refGE}[1]{\ensuremath{\stackrel{(\ref{#1})}{\geq}}}

\newcommand{\refPE}[1]{\ensuremath{\stackrel{(\ref{#1})}{\preceq}}}
\newcommand{\refSE}[1]{\ensuremath{\stackrel{(\ref{#1})}{\succeq}}}

\def\Det{{\rm Det}}
\def\Tr{{\rm Tr}}
\def\DFP{{\rm DFP}}
\def\BFGS{{\rm BFGS}}

\def\Broyd{{\rm Broyd}}

\renewcommand\arraystretch{1.5}
\def\ba{\begin{array}}
\def\ea{\end{array}}
\def\beann{\begin{eqnarray*}}
\def\eeann{\end{eqnarray*}}
\def\bea{\begin{eqnarray}}
\def\eea{\end{eqnarray}}
\def\bal{\renewcommand\arraystretch{2}\begin{tabular}{|l|}}
\def\eal{\end{tabular}}

\def\BT{\begin{theorem}}
\def\ET{\end{theorem}}
\def\BL{\begin{lemma}}
\def\EL{\end{lemma}}
\def\BC{\begin{corollary}}
\def\EC{\end{corollary}}
\def\BE{\begin{example}}
\def\EE{\end{example}}
\def\BD{\begin{definition}}
\def\ED{\end{definition}}
\def\BR{\begin{remark}}
\def\ER{\end{remark}}
\def\BAS{\begin{assumption}}
\def\EAS{\end{assumption}}
\def\BI{\begin{itemize}}
\def\EI{\end{itemize}}
\def\BCA{\begin{cases}}
\def\ECA{\end{cases}}

\def\BMP{\begin{minipage}{9.5cm}}
\def\EMP{\end{minipage}}
\def\MPT{\begin{minipage}{11.5cm}}
\def\EPT{\end{minipage}}

\def\la{\langle}
\def\ra{\rangle}

\def\qedEA{\tag*{$\qed$}}

\begin{document}

\title{New Results on Superlinear Convergence of Classical
Quasi-Newton Methods\thanks{Communicated by Russell Luke.}}

\author{Anton Rodomanov   \and  Yurii Nesterov}

\institute{Anton Rodomanov \at ICTEAM,
Catholic University of Louvain, Louvain-la-Neuve, Belgium \\
anton.rodomanov@uclouvain.be
\and Yurii Nesterov \at CORE, Catholic University of
Louvain, Louvain-la-Neuve, Belgium \\
yurii.nesterov@uclouvain.be}

\date{%
   Received: 14 July 2020 / Accepted: 25 December 2020
   / Published online: 9 January 2021 \\
   \textcopyright{} The Author(s) 2021
}
%The correct dates will be entered by the editor.

\def\subclassname{{\bfseries Mathematics Subject Classification}\enspace}
\headerboxheight=70pt
\def\makeheadbox{\noindent\small
  Journal of Optimization Theory and Applications (2021) 188:744--769 \\
  \url{https://doi.org/10.1007/s10957-020-01805-8}
}

\maketitle

\begin{abstract}
We present a new theoretical analysis of local superlinear
convergence of classical quasi-Newton methods from the
convex Broyden class. As a result, we obtain a significant
improvement in the currently known estimates of the
convergence rates for these methods. In particular, we show
that the corresponding rate of the
Broyden--Fletcher--Goldfarb--Shanno method depends
only on the product of the dimensionality of the problem and
the \emph{logarithm} of its condition number.
\end{abstract}


\section{Introduction}

We study local superlinear convergence of classical
quasi-Newton methods for smooth unconstrained optimization.
These algorithms can be seen as an approximation of the
standard Newton method, in which the exact Hessian is
replaced by some operator, which is updated in iterations by
using the gradients of the objective function. The two most
famous examples of quasi-Newton algorithms are the
Davidon--Fletcher--Powell (DFP)
\cite{Davidon1959,FletcherPowell1963} and the
Broyden--Fletcher--Goldfarb--Shanno (BFGS)
\cite{Broyden1970p1,Broyden1970p2,Fletcher1970,Goldfarb1970,Shanno1970}
methods, which together belong to the Broyden family
\cite{Broyden1967} of quasi-Newton algorithms. For an
introduction into the topic, see \cite{DennisMore1977} and
\cite[Chapter~6]{NocedalWright2006}. See also
\cite{LewisOverton2013} for the discussion of quasi-Newton
algorithms in the context of nonsmooth optimization.

The superlinear convergence of quasi-Newton methods was
established as early as in 1970s, firstly by Powell
\cite{Powell1971} and Dixon \cite{Dixon1972p1,Dixon1972p2}
for the methods with exact line search, and then by
Broyden, Dennis and Mor{\'e} \cite{BroydenDennisMore1973}
and Dennis and Mor{\'e} \cite{DennisMore1974} for the
methods without line search. The latter two approaches
have been extended onto more general methods under various
settings (see, e.g.,
\cite{Stachurski1981,GriewankToint1982,EngelsMartinez1991,ByrdLiuNocedal1992,YabeYamaki1996,WeiYuYuanLian2004,YabeOgasawaraYoshino2007,MokhtariEisenRibeiro2018,GaoGoldfarb2019}).

However, explicit \emph{rates} of superlinear convergence
for quasi-Newton algorithms were obtained only
recently. The first results were presented in
\cite{RodomanovNesterov2020a} for the \emph{greedy}
quasi-Newton methods. After that, in
\cite{RodomanovNesterov2020b}, the \emph{classical}
quasi-Newton methods were considered, for which the authors
established certain superlinear convergence rates, depending
on the problem dimension and its condition number. The
analysis was based on the trace potential function, which
was then augmented by the logarithm of determinant of the
\emph{inverse} Hessian approximation to extend the proof
onto the general nonlinear case.

In this paper, we further improve the results of
\cite{RodomanovNesterov2020b}. For the classical
quasi-Newton methods, we obtain new convergence rate
estimates, which have better dependency on the condition
number of the problem. In particular, we show that the
superlinear convergence rate of BFGS depends on the
condition number only through the \emph{logarithm}. As
compared to the previous work, the main difference in the
analysis is the choice of the potential function: now the
main part is formed by the logarithm of determinant of
Hessian approximation, which is then augmented by the trace
of \emph{inverse} Hessian approximation.

It is worth noting that recently, in \cite{JinMokhtari2020},
another analysis of local superlinear convergence of the
classical DFP and BFGS methods was presented with the
resulting rate, which is independent of the dimensionality
of the problem and its condition number. However, to obtain
such a rate, the authors had to make an additional
assumption that the methods start from a sufficiently good
initial Hessian approximation. Without this assumption, to
our knowledge, their proof technique, based on the
Frobenius-norm potential function, leads only to the rates,
which are weaker than those in
\cite{RodomanovNesterov2020b}.

This paper is organized as follows. In Sect.~\ref{sec-not},
we introduce our notation. In Sect.~\ref{sec-broyd}, we
study the convex Broyden class of quasi-Newton updates for
approximating a self-adjoint positive definite operator. In
Sect.~\ref{sec-quad}, we analyze the rate of convergence of
the classical quasi-Newton methods from the convex Broyden
class as applied to minimizing a quadratic function. On this
simple example, where the Hessian is constant, we illustrate
the main ideas of our analysis. In Sect.~\ref{sec-gen}, we
consider the general unconstrained optimization problem.
Finally, in Sect.~\ref{sec-disc}, we discuss why the new
superlinear convergence rates, obtained in this paper, are
better than the previously known ones.

\section{Notation}\label{sec-not}

In what follows, $\E$ denotes an $n$-dimensional
real vector space. Its dual space, composed of all linear
functionals on $\E$, is denoted by $\E^*$. The value of a
linear function $s \in \E^*$, evaluated at a point $x \in
\E$, is denoted by $\la s, x \ra$.

For a smooth function $f : \E \to \R$, we denote by $\nabla
f(x)$ and $\nabla^2 f(x)$ its gradient and Hessian
respectively, evaluated at a point $x \in \E$. Note that
$\nabla f(x) \in \E^*$, and $\nabla^2 f(x)$ is a
self-adjoint linear operator from $\E$ to $\E^*$.

The partial ordering of self-adjoint linear operators is
defined in the standard way. We write $A_1 \preceq A_2$ for
$A_1, A_2 : \E \to \E^*$, if $\la (A_2 - A_1) x, x
\ra \geq 0$ for all $x \in \E$, and $H_1 \preceq H_2$ for
$H_1, H_2 : \E^* \to \E$, if $\la s, (H_2 - H_1) s
\ra \geq 0$ for all $s \in \E^*$.

Any self-adjoint positive definite linear operator $A : \E
\to \E^*$ induces in the spaces $\E$ and $\E^*$ the
following pair of conjugate Euclidean norms:
\beq\label{def-norms}
\ba{rclrcl}
\| h \|_A &\Def& \la A h, h \ra^{1/2}, \quad h \in \E,
\qquad\quad
\| s \|_A^* &\Def& \la s, A^{-1} s \ra^{1/2}, \quad s \in
\E^*.
\ea
\eeq
When $A = \nabla^2 f(x)$, where $f : \E \to \R$ is a smooth
function with positive definite Hessian, and $x \in \E$, we
prefer to use notation $\| \cdot \|_x$ and $\| \cdot
\|_x^*$, provided that there is no ambiguity with the
reference function $f$.

Sometimes, in the formulas, involving products of linear
operators, it is convenient to treat $x \in \E$ as a
linear operator from $\R$ to $\E$, defined by $x \alpha =
\alpha x$, and $x^*$ as a linear operator from $\E^*$ to
$\R$, defined by $x^* s = \la s, x \ra$. Likewise, any $s
\in \E^*$ can be treated as a linear operator from $\R$ to
$\E^*$, defined by $s \alpha = \alpha s$, and $s^*$ as a
linear operator from $\E$ to $\R$, defined by $s^* x = \la
s, x \ra$. In this case, $x x^*$ and $s s^*$ are rank-one
self-adjoint linear operators from~$\E^*$ to $\E$ and
from~$\E^*$ to $\E$ respectively, acting as follows: $(x x^*) s
= \la s, x \ra x$ and $(s s^*) x = \la s, x \ra s$ for $x
\in \E$ and $s \in \E^*$.

Given two self-adjoint linear operators $A : \E \to \E^*$
and $H : \E^* \to \E$, we define the trace and the
determinant of $A$ with respect to $H$ as follows: $\la H, A
\ra \Def \Tr(H A)$, and $\Det(H, A) \Def \Det(H A)$. Note
that $H A$ is a linear operator from $\E$ to itself, and
hence its trace and determinant are well-defined by the
eigenvalues (they coincide with the trace and determinant of
the matrix representation of $H A$ with respect to an
arbitrary chosen basis in the space~$\E$, and the result is
independent of the particular choice of the basis). In
particular, if $H$ is positive definite, then $\la H, A \ra$
and $\Det(H, A)$ are respectively the sum and the product of
the eigenvalues of $A$ relative to $H^{-1}$. Observe that
$\la \cdot, \cdot \ra$ is a bilinear form, and for any $x
\in \E$, we have $\la A x, x \ra = \la x x^*, A \ra$. When
$A$ is invertible, we also have $\la A^{-1}, A \ra = n$ and
$\Det(A^{-1}, \delta A) = \delta^n$ for any $\delta \in \R$.
Also recall the following multiplicative formula for the
determinant: $\Det(H, A) = \Det(H, G) \cdot \Det(G^{-1},
A)$, which is valid for any invertible linear operator $G :
\E \to \E^*$. If the operator $H$ is positive
semidefinite, and $A_1 \preceq A_2$ for some self-adjoint
linear operators $A_1, A_2 : \E \to \E^*$, then $\la H, A_1
\ra \leq \la H, A_2 \ra$ and $\Det(H, A_1) \leq \Det
(H, A_2)$. Similarly, if $A$ is positive semidefinite and
$H_1 \preceq H_2$ for some self-adjoint linear operators
$H_1, H_2 : \E^* \to \E$, then $\la H_1, A \ra \leq \la H_2,
A \ra$ and $\Det(H_1, A) \leq \Det(H_2, A)$.

\section{Convex Broyden Class}\label{sec-broyd}

Let $A$ and $G$ be two self-adjoint positive definite linear
operators from $\E$ to~$\E^*$, where $A$ is the target
operator, which we want to approximate, and $G$ is its
current approximation. The \emph{Broyden class} of
quasi-Newton updates of~$G$ with respect to $A$ along a
direction $u \in \E \setminus \{0\}$ is the following family
of updating formulas, parameterized by a scalar $\tau \in
\R$:
\beq\label{def-broyd}
\ba{rcl}
\Broyd_{\tau}(A, G, u) &=& \phi_{\tau} \left[ G -
\frac{A u u^* G + G u u^* A}{\la A u, u \ra} + \left(
\frac{\la G u, u \ra}{\la A u, u \ra} + 1 \right) \frac{A u
u^* A}{\la A u, u \ra} \right] \\
&& + \, (1 - \phi_{\tau}) \left[ G - \frac{G u u^*
G}{\la G u, u \ra} + \frac{A u u^* A}{\la A u, u \ra}
\right],
\ea
\eeq
where
\beq\label{def-phi}
\ba{rcl}
\phi_{\tau} \;\Def\; \phi_{\tau}(A, G, u) &\Def& \frac{\tau
\frac{\la A u, u \ra}
{\la A G^{-1} A u, u \ra}}{\tau \frac{\la A u, u \ra}{\la A
G^{-1} A u, u \ra} + (1 - \tau) \frac{\la G u, u \ra}{\la A
u, u \ra}}.
\ea
\eeq
If the denominator in \eqref{def-phi} is zero, we left
both $\phi_{\tau}$ and $\Broyd_{\tau}(A, G, u)$ undefined.
For the sake of convenience, we also set $\Broyd_{\tau}
(A, G, u) = G$ for $u = 0$.

In this paper, we are interested in the \emph{convex}
Broyden class, which is described by the values of $\tau \in
[0, 1]$. Note that for all such $\tau$ the denominator in
\eqref{def-phi} is always positive for any $u \neq
0$, so both $\phi_{\tau}$ and $\Broyd_{\tau}(A, G, u)$ are
well-defined; moreover, $\phi_{\tau} \in [0, 1]$. For $\tau
= 1$, we have $\phi_{\tau} = 1$, and \eqref{def-broyd}
becomes the DFP update; for $\tau = 0$, we have
$\phi_{\tau} = 0$, and \eqref{def-broyd} becomes the BFGS
update.

\BR
Usually the Broyden class is defined directly in terms of
the parameter $\phi$. However, in the context of this paper,
it is more convenient to work with $\tau$ instead of $\phi$.
As can be seen from \eqref{inv-broyd}, $\tau$ is exactly the
weight of the DFP component in the updating formula for the
inverse operator.
\ER

A basic property of an update from the convex Broyden
class is that it preserves the bounds on the eigenvalues
with respect to the target operator.

\BL[see~{\cite[Lemma~2.1]{RodomanovNesterov2020b}}]
\label{lm-eigs}
If $\frac{1}{\xi} A \preceq G \preceq \eta A$ for some
$\xi, \eta \geq 1$, then, for any $u \in \E$, and any $\tau
\in [0, 1]$, we have $\frac{1}{\xi} A \preceq
\Broyd_{\tau}(A, G, u) \preceq \eta A$.
\EL

Consider the measure of closeness of $G$ to $A$ along
direction $u \in \E \setminus \{0\}$:
\beq\label{def-nu}
\ba{rcl}
\nu(A, G, u) &\Def& \frac{\la (G - A) G^{-1} (G - A) u, u
\ra^{1/2}}{\la A u, u \ra^{1/2}} \;\refEQ{def-norms}\;
\frac{\| (G - A) u \|_G^*}{\| u \|_A}.
\ea
\eeq
Let us present two potential functions, whose improvement
after one update from the convex Broyden class can be
bounded from below by a certain non-negative monotonically
increasing function of $\nu$, vanishing at zero.

First, consider the \emph{log-det barrier}
\beq\label{def-V}
\ba{rcl}
V(A, G) &=& \ln \Det(A^{-1}, G).
\ea
\eeq
It will be useful when $A \preceq G$. Note that in this case
$V(A, G) \geq 0$.

\BL\label{lm-prog-V}
Let $A, G : \E \to \E^*$ be self-adjoint positive definite
linear operators, $A \preceq G \preceq \eta A$ for
some $\eta \geq 1$. Then, for any $\tau \in [0, 1]$ and $u
\in \E \setminus \{0\}$:
$$
\ba{rcl}
V(A, G) - V(A, \Broyd_{\tau}(A, G, u)) &\geq& \ln\left( 1 +
( \tau \frac{1}{\eta} + 1 - \tau) \nu^2(A, G, u) \right).
\ea
$$
\EL

\begin{proof}
Indeed, denoting $G_+ \Def \Broyd_{\tau}(A, G, u)$, we
obtain
\beq\label{prog-V-prel}
\ba{rcl}
&&\hspace{-6em} V(A, G) - V(A, G_+) \;\refEQ{def-V}\; 
\ln\Det(G_+^{-1}, G) \\
&\refEQ{det-broyd}& \ln\left( \tau \frac{\la A u,
u \ra}{\la A G^{-1} A u, u \ra} + (1 - \tau) \frac{\la G u,
u \ra}{\la A u, u \ra} \right) \\
&=& \ln\left( 1 + \tau \frac{\la A (A^{-1} - G^{-1}) A u, u
\ra}{\la A G^{-1} A u, u \ra} + (1 - \tau) \frac{\la (G - A)
u, u \ra}{\la A u, u \ra} \right).
\ea
\eeq
Since\footnote{This is obvious when $G - A$ is
non-degenerate. The general case follows by continuity.} $0
\preceq G - A \preceq (1 - \frac{1}{\eta}) G$, we have
\beq\label{prog-V-aux}
\ba{rcl}
(G - A) G^{-1} (G - A) &\preceq& \left( 1 - \frac{1}{\eta}
\right) (G - A) \;\preceq\; \frac{1}{1 + \frac{1}{\eta}} (G
- A) \;\preceq\; G - A.
\ea
\eeq
Therefore, denoting $\nu \Def \nu(A, G, u)$, we can write
that
$$
\ba{rcl}
\frac{\la (G - A) u, u \ra}{\la A u, u \ra} &\refGE
{prog-V-aux}& \frac{\la (G - A) G^{-1} (G - A) u, u
\ra}{\la A u, u \ra} \;\refEQ{def-nu}\; \nu^2,
\ea
$$
and, since $A (A^{-1} - G^{-1}) A = G - A - (G - A) G^{-1} 
(G - A)$, that
$$
\ba{rcl}
\frac{\la A (A^{-1} - G^{-1}) A u, u \ra}{\la A G^{-1} A u,
u \ra} &=& \frac{\la (G - A - (G - A) G^{-1}
(G - A)) u, u \ra}{\la A G^{-1} A u, u \ra} \;\refGE
{prog-V-aux}\; \frac{1}{\eta} \frac{\la (G - A) G^{-1} (G -
A) u, u \ra}{\la A G^{-1} A u, u \ra} \\
&\geq& \frac{1}{\eta} \frac{\la (G - A)
G^{-1} (G - A) u, u \ra}{\la A u, u \ra} \;\refEQ{def-nu}\;
\frac{1}{\eta} \nu^2.
\ea
$$
Substituting the above two inequalities into
\eqref{prog-V-prel}, we obtain the claim.\qed
\end{proof}

Now consider another potential function, the \emph{augmented
log-det barrier}:
\beq\label{def-psi}
\ba{rcl}
\psi(G, A) &\Def& \ln\Det(A^{-1}, G) - \la G^{-1}, G - A
\ra.
\ea
\eeq
As compared to the log-det barrier, this potential
function is more universal since it works even if the
condition $A \preceq G$ is violated. Note that the
augmented log-det barrier is in fact the Bregman
divergence, generated by the strictly convex function
$d(A) \Def -\ln\Det(B^{-1}, A)$, defined on the set of
self-adjoint positive definite linear operators from $\E$
to $\E^*$, where $B : \E \to \E^*$ is an arbitrary fixed
self-adjoint positive definite linear operator. Indeed,
\beq\label{psi-breg}
\ba{rcl}
\psi(G, A) &=& -\ln\Det(B^{-1}, A) +
\ln\Det(B^{-1}, G) - \la -G^{-1}, A - G \ra \\
&=& d(A) - d(G) - \la \nabla d(G), A - G \ra \;\geq\; 0.
\ea
\eeq

\BR
The idea of combining the trace with the logarithm of
determinant to form a potential function for the analysis of
quasi-Newton methods can be traced back to
\cite{ByrdNocedal1989}. Note also that in
\cite{RodomanovNesterov2020b}, the authors studied the
evolution of $\psi(A, G)$, i.e. the Bregman divergence was
centered at $A$ instead of $G$.
\ER

\BL\label{lm-aux-ineq}
For any real $\alpha \geq \beta > 0$, we have $\alpha +
\frac{1}{\beta} - 1 \geq 1$, and
\beq\label{aux-ineq}
\ba{rclrcl}
\alpha - \ln \beta - 1 &\geq& \frac{\sqrt{3}}{2 + \sqrt{3}}
\ln\left( \alpha + \frac{1}{\beta} - 1 \right) &\geq&
\frac{6}{13} \ln\left( \alpha + \frac{1}{\beta} - 1 \right).
\ea
\eeq
\EL

\begin{proof}
We only need to prove the first inequality in
\eqref{aux-ineq} since the second one follows from it and
the fact that $\frac{\sqrt{3} + 2} {\sqrt{3}} = 1 +
\frac{2}{\sqrt{3}} \leq 1 + \frac{7}{6} =
\frac{13}{6}$ (since $2 \leq \frac{7}{2 \sqrt{3}}$).

Let $\beta > 0$ be fixed, and let $\zeta_1 : (1 - \frac{1}
{\beta}, +\infty) \to \R$ be the function, defined by
$\zeta_1(\alpha) \Def \alpha - \frac{\sqrt{3}}{2 +
\sqrt{3}} \ln\left( \alpha + \frac{1}{\beta} - 1 \right)$.
Note that the domain of $\zeta_1$ includes the point
$\alpha = \beta$ since $\beta \geq 2 - \frac{1}{\beta} > 1 -
\frac{1}{\beta}$. Let us show that $\zeta_1$ increases on
the interval $[\beta, +\infty)$. Indeed, for any $\alpha
\geq \beta$, we have
$$
\ba{rcl}
\zeta_1'(\alpha) &=& 1 - \frac{\sqrt{3}}{2 + \sqrt{3}} 
\frac{1}{\alpha + \frac{1}{\beta} - 1} \;>\;
1 - \frac{1}{\alpha + \frac{1}{\beta} - 1} \;=\; 
\frac{\alpha + \frac{1}{\beta} - 2}{\alpha + \frac{1}{\beta}
- 1} \;\geq\; \frac{\beta + \frac{1}{\beta} - 2}{\alpha +
\frac{1}{\beta} - 1} \;\geq\; 0.
\ea
$$
Thus, it is sufficient to prove \eqref{aux-ineq} only in the
case when $\alpha = \beta$. Equivalently, we need to
show that the function $\zeta_2 : (0, +\infty) \to \R$,
defined by the formula $\zeta_2 (\alpha) \Def \alpha - \ln
\alpha - 1 -
\frac{
\sqrt{3}}{2 + \sqrt{3}} \ln\left( \alpha + \frac{1}{\alpha}
- 1 \right)$, is non-negative. Differentiating, we find
that, for all $\alpha > 0$, we have
$$
\ba{rcl}
\zeta_2'(\alpha) &=& 1 - \frac{1}{\alpha} -
\frac{\sqrt{3}}{2 + \sqrt{3}} \frac{1 - \frac{1}{\alpha^2}}
{\alpha + \frac{1}{\alpha} - 1} \;=\; \left( 1 - \frac{1}
{\alpha} \right) \left( 1 - \frac{\sqrt{3}}{2 + \sqrt{3}}
\frac{1 + \frac{1}{\alpha}}{\alpha + \frac{1}{\alpha} - 1}
\right) \\
&=& \left( 1 - \frac{1}{\alpha} \right) \frac{\alpha +
\frac{1}{\alpha} - 1 - (2 \sqrt{3} - 3) (1 + \frac{1}
{\alpha})}{\alpha + \frac{1}{\alpha} - 1} \;=\; \left( 1 -
\frac{1}{\alpha} \right) \frac{\alpha - 2 (\sqrt{3} - 1) +
(\sqrt{3} - 1)^2 \frac{1}{\alpha}}{1 + \frac{1}{\alpha} -
1} \\
&=& \left( 1 - \frac{1}{\alpha} \right) \frac{(
\sqrt{\alpha} - (\sqrt{3} - 1) \frac{1}{\sqrt{\alpha}})^2}
{\alpha + \frac{1}{\alpha} - 1}.
\ea
$$
Hence, $\zeta_2'(\alpha) \leq 0$ for $0 < \alpha \leq 1$,
and $\zeta_2'(\alpha) \geq 0$ for $\alpha \geq 1$. Thus,
the minimum of $\zeta_2$ is attained at $\alpha = 1$.
Consequently, $\zeta_2(\alpha) \geq \zeta_2(1) = 0$ for all
$\alpha > 0$.\qed
\end{proof}

It turns out that, up to some constants, the improvement in
the augmented log-det barrier can be bounded from below by
exactly the same logarithmic function of $\nu$, which was
used for the simple log-det barrier.

\BL\label{lm-prog-psi}
Let $A, G : \E \to \E^*$ be self-adjoint positive definite
linear operators, $\frac{1}{\xi} A \preceq G
\preceq \eta A$ for some $\xi, \eta \geq 1$. Then, for any
$\tau \in [0, 1]$ and $u \in \E \setminus \{0\}$:
$$
\ba{rcl}
\psi(G, A) - \psi(\Broyd_{\tau}(A, G, u), A) &\geq&
\frac{6}{13} \ln\left(1 + ( \tau \frac{1}
{\xi \eta} + 1 - \tau ) \nu^2(A, G, u) \right).
\ea
$$
\EL

\begin{proof}
Indeed, denoting $G_+ \Def \Broyd_{\tau}(A, G, u)$, we
obtain
$$
\ba{rcl}
\la G^{-1} - G_+^{-1}, A \ra &\refEQ{inv-broyd}& \tau \left[
\frac{\la A G^{-1} A G^{-1} A u, u \ra} {\la A G^{-1} A u, u
\ra} - 1 \right] + (1 - \tau) \left[ \frac{\la A G^{-1} A u,
u \ra}{\la A u, u \ra} - 1 \right],
\ea
$$
and
$$
\ba{rcl}
\Det(G_+^{-1}, G) &\refEQ{det-broyd}& \tau \frac{\la A u, u
\ra}{\la A G^{-1} A u, u \ra} + (1 - \tau) \frac{\la A u, u
\ra}{\la G u, u \ra}.
\ea
$$
Thus,
\beq\label{prog-psi-prel}
\ba{rcl}
&&\hspace{-1em} \psi(G, A) - \psi(G_+, A) \;\refEQ
{def-psi}\; \la G^{-1}
- G_+^{-1}, A \ra + \ln\Det (G_+^{-1}, G) \\
&&\hspace{1em} \;=\; \tau \alpha_1 + (1 - \tau) \alpha_0 +
\ln(\tau
\beta_1^
{-1} + (1 - \tau) \beta_0^{-1} ) - 1 \;=\; \alpha - \ln
\beta - 1,
\ea
\eeq
where we denote $\alpha_1 \Def \frac{\la A G^{-1} A G^
{-1} A u, u
\ra}{\la A G^{-1} A u, u \ra}$, $\beta_1 \Def \frac{\la A
G^{-1} A u, u \ra}{\la A u, u \ra}$, $\alpha_0 \Def 
\frac{\la A G^{-1} A u, u \ra}{\la A u, u \ra}$, $\beta_0
\Def \frac{\la A u, u \ra}{\la G u, u \ra}$, $\alpha \Def
\tau \alpha_1 + (1 - \tau) \alpha_0$, $\beta \Def ( \tau
\beta_1^{-1} + (1 - \tau) \beta_0^{-1} )^{-1}$. Note that
$\alpha_1 \geq \beta_1$ and $\alpha_0 \geq
\beta_0$ by the Cauchy-Schwartz inequality. At the
same time, $\tau \beta_1 + (1 - \tau) \beta_2 \geq \beta$ by
the convexity of the inverse function $t \mapsto t^{-1}$.
Hence, we can apply Lemma~\ref{lm-aux-ineq} to estimate
\eqref{prog-psi-prel} from below. Note that
\[
\ba[b]{rcl}
\alpha + \frac{1}{\beta} - 1 &=& \tau
\frac{\la (A + A G^{-1} A G^{-1} A) u, u \ra}{\la A G^ {-1}
A u, u \ra} + (1 - \tau) \frac{\la (G + A) u, u \ra} {\la A
u, u \ra} - 1 \\
&=& 1 + \tau \frac{\la (G - A) G^{-1} A G^{-1} (G - A) \ra}
{\la A G^{-1} A u, u \ra} + (1 - \tau) \frac{\la (G - A) G^
{-1} (G - A) u, u \ra}{\la A u, u \ra} \\
&\geq& 1 + (\tau \frac{1}{\xi \eta} + 1 -
\tau) \frac{\la (G - A) G^{-1} (G - A) u, u \ra}{\la A u, u
\ra} \\
&\refEQ{def-nu}& 1 + (\tau \frac{1}{\xi \eta} + 1 - \tau)
\nu^2(A, G, u).
\ea\qedEA
\]
\end{proof}

The measure $\nu(A, G, u)$, defined in \eqref{def-nu}, is
the ratio of the norm of $(G - A) u$, measured with
respect to $G$, and the norm of $u$, measured with respect
to $A$. It is important that we can change the
corresponding metrics to $G_+$ and $G$ respectively by
paying only with the minimal eigenvalue of $G$ relative to
$A$.

\BL\label{lm-ch-met}
Let $A, G : \E \to \E^*$ be self-adjoint positive definite
linear operators such that $\frac{1}{\xi} A \preceq G$ for
some $\xi > 0$. Then, for any $\tau \in [0, 1]$, any $u
\in \E \setminus \{0\}$, and $G_+ \Def \Broyd_{\tau}(A, G,
u)$, we have
$$
\ba{rcl}
\nu^2(A, G, u) &\geq& \frac{1}{1 + \xi} \frac{\la (G -
A) G_+^{-1} (G - A) u, u \ra}{\la G u, u \ra}.
\ea
$$
\EL

\begin{proof}
From \eqref{inv-broyd}, it is easy to see that $G_+^{-1} A u
= u$. Hence,
\beq\label{ch-met-aux1}
\ba{rcl}
\frac{\la (G - A) G_+^{-1} (G - A) u, u
\ra}{\la G u, u \ra} &=& \frac{\la G G_+^{-1} G u, u
\ra}{\la G u, u \ra} +
\frac{\la A u, G_+^{-1} A u \ra}{\la G u, u \ra} - 2
\frac{\la G u, G_+^{-1} A u \ra}{\la G u, u \ra} \\
&=& \frac{\la G G_+^{-1} G u, u \ra}{\la G u, u
\ra} + \frac{\la A u, u \ra}{\la G u, u \ra} - 2.
\ea
\eeq
Since $1 - t \leq \frac{1}{t} - 1$ for all $t > 0$, we
further have
\beq\label{G-Gp-G}
\ba{rcl}
\frac{\la G G_+^{-1} G u, u \ra}{\la G u, u \ra}
&\refEQ{inv-broyd}& \tau \left[ 1 - \frac{\la A u, u \ra^2}
{\la G u, u \ra \la A G^{-1} A u, u \ra} + \frac{\la G u, u
\ra}{\la A u, u \ra} \right] \\
&& + \, (1 - \tau) \left[ \left( \frac{\la A G^{-1} A u, u
\ra}
{\la A u, u \ra} + 1 \right) \frac{\la G u, u \ra}{\la A u,
u \ra} - 1 \right] \\
&\leq& \left( \frac{\la A G^{-1} A u, u \ra}{\la A u, u \ra}
+ 1 \right) \frac{\la G u, u \ra}{\la A u, u \ra} - 1.
\ea
\eeq
Denote $\nu \Def \nu(A, G, u)$. Then,
\beq\label{ch-met-nu2}
\ba{rcl}
\nu^2 &\refEQ{def-nu}& \frac{\la (G - A) G^{-1} (G - A) u,
u \ra}{\la A u, u \ra} \;=\; \frac{\la G u, u \ra}{\la A u,
u \ra} + \frac{\la A G^{-1} A u, u \ra}{\la A u, u \ra} - 2.
\ea
\eeq
Consequently,
\beq\label{ch-met-aux2}
\ba{rcl}
\hspace{-1em} (1 + \xi) \nu^2 &\geq& \left( 
\frac{\la A G^{-1} A u, u \ra}{\la A u, u \ra} + 1 \right)
\nu^2 \\ &\refEQ{ch-met-nu2}& \left( \frac{\la A G^{-1} A u,
u \ra}{\la A u, u \ra} + 1 \right) \frac{\la G u, u \ra}{\la
A u, u \ra} + \frac{\la A G^{-1} A u, u \ra^2}{\la A u, u
\ra^2} - \frac{\la A G^{-1} A u, u \ra}{\la A u, u \ra} - 2
\\
&\refGE{G-Gp-G}& \frac{\la G G_+^{-1} G u, u \ra}{\la A u,
u \ra} + \frac{\la A G^{-1} A u, u \ra^2}{\la A u, u \ra}
- \frac{\la A G^{-1} A u, u \ra}{\la A u, u \ra} - 1.
\ea
\eeq
Thus,
$$
\ba{rcl}
(1 + \xi) \nu^2 - \frac{\la (G - A) G_+^{-1} (G - A) u, u
\ra}{\la G u, u \ra} &\refEQ{ch-met-aux1}& (1 + \xi) \nu^2 -
\frac{\la G G_+^{-1} G u, u \ra}{G u, u \ra} - \frac{\la A
u, u \ra}{\la G u, u \ra} + 2 \\
&\refGE{ch-met-aux2}& \frac{\la A G^{-1} A u, u \ra^2}{\la
A u, u \ra^2} - \frac{\la A G^{-1} A u, u \ra}{\la A u, u
\ra} - \frac{\la A u, u \ra}{\la G u, u \ra} + 1 \\
&\geq& \frac{\la A G^{-1} A u, u \ra^2}{\la A u, u \ra^2} -
2 \frac{\la A G^{-1} A u, u \ra}{\la A u, u \ra} + 1
\;\geq\; 0,
\ea
$$
where we have used the Cauchy--Schwartz inequality $
\frac{\la A u, u \ra}{\la G u, u \ra} \leq \frac{\la A G^
{-1} A u, u \ra}{\la A u, u \ra}$.\qed
\end{proof}

\section{Unconstrained Quadratic Minimization}
\label{sec-quad}

Let us study the convergence properties of the classical
quasi-Newton methods from the convex Broyden class, as
applied to minimizing the quadratic function
\beq\label{func-quad}
\ba{rcl}
f(x) &\Def& \frac{1}{2} \la A x, x \ra - \la b, x \ra,
\ea
\eeq
where $A : \E \to \E^*$ is a self-adjoint positive definite
linear operator, and $b \in \E^*$.

Let $B : \E \to \E^*$ be a fixed self-adjoint positive
definite linear operator, and let $\mu, L > 0$ be such that
\beq\label{quad-mu-L}
\ba{rclrcl}
\mu B &\preceq& A &\preceq& L B.
\ea
\eeq
Thus, $\mu$ is the \emph{strong convexity} parameter of $f$,
and $L$ is the constant of \emph{Lipschitz continuity} of
the gradient of $f$, both measured relative to $B$.

Consider the following standard quasi-Newton process for
minimizing \eqref{func-quad}:
\beq\label{qn-quad}
\bal \hline
\textbf{Initialization:} Choose $x_0 \in \E$. Set $G_0 = L
B$. \\ \hline \textbf{For $k \geq 0$ iterate:} \\
1. Update $x_{k+1} = x_k - G_k^{-1} \nabla f(x_k)$. \\
2. Set $u_k = x_{k+1} - x_k$ and choose $\tau_k \in [0, 1]$.
\\
3. Compute $G_{k+1} = \Broyd_{\tau_k}(A, G_k, u_k)$.
\\
\hline
\eal
\eeq
For measuring its rate of convergence, we use the norm of
the gradient, taken with respect to the Hessian:
$$
\ba{rcl}
\lambda_k &\Def& \| \nabla f(x_k) \|_A^* \;\refEQ
{def-norms}\; \la \nabla f(x_k), A^{-1} \nabla f(x_k) \ra^
{1/2}.
\ea
$$

It is known that the process \eqref{qn-quad} has at least
a linear convergence rate of the standard gradient method:
\BT[see {\cite[Theorem~3.1]{RodomanovNesterov2020b}}]
\label{th-lin-quad}
In scheme~\eqref{qn-quad}, for all $k \geq 0$:
\beq\label{quad-lin}
\ba{rclrclrcl}
A &\preceq& G_k &\preceq& \frac{L}{\mu} A,
\qquad
\lambda_k &\leq& \left( 1 - \frac{\mu}{L} \right)^k \,
\lambda_0.
\ea
\eeq
\ET

Let us establish the superlinear convergence. According to
\eqref{quad-lin}, for the quadratic function, we have $A
\preceq G_k$ for all $k \geq 0$. Therefore, in our
analysis, we can use both potential functions: the log-det
barrier and the augmented log-det barrier. Let us consider
both options. We start with the first one.

\BT\label{th-sup-quad}
In scheme~\eqref{qn-quad}, for all $k \geq 1$, we have
\beq\label{sup-quad}
\ba{rcl}
\lambda_k &\leq& \left[\frac{2}{\prod_{i=0}^{k-1} (\tau_i
\frac{\mu}{L} + 1 - \tau_i)^{1/k}} \left (e^{\frac{n}{k} \ln
\frac{L}{\mu}} - 1 \right) \right]^{k/2}
\sqrt{\frac{L}{\mu}} \cdot \lambda_0.
\ea
\eeq
\ET

\begin{proof}
Without loss of generality, we can assume that $u_i \neq 0$
for all $0 \leq i \leq k$. Denote $V_i \Def V(A,
G_i)$, $\nu_i \Def \nu(A, G_i, u_i)$, $p_i \Def \tau_i
\frac{\mu}{L} + 1 - \tau_i$, $g_i \Def \| \nabla f(x_i) \|_
{G_i}^*$ for any $0 \leq i \leq k$. By
Lemma~\ref{lm-prog-V} and \eqref{quad-lin}, for all $0 \leq i
\leq k - 1$, we have $\ln( 1 + p_i \nu_i^2 ) \leq V_i - V_
{i+1}$. Summing up, we obtain
\beq\label{sup-prel-quad}
\ba{rcl}
&&\hspace{-5em}\sum\limits_{i=0}^{k-1} \ln(1 + p_k \nu_k^2)
\;\leq\; V_0 - V_k \;\refLE{quad-lin}\; V_0 \;\refEQ
{qn-quad}\; V(A, L B) \\
&\refEQ{def-V}& \ln \Det(A^{-1}, L B) \;\refLE{quad-mu-L}\;
\ln \Det(\frac{1}{\mu} B^{-1}, L B)
\;=\; n \ln \frac{L}{\mu}.
\ea
\eeq
Hence, by the convexity of function $t \mapsto \ln(1 +
e^t)$, we get
\beq\label{sup-prel-quad-2}
\ba{rcl}
\frac{n}{k} \ln \frac{L}{\mu} &\refGE{sup-prel-quad}&
\frac{1}{k} \sum\limits_{i=0}^{k-1} \ln(1 + p_i \nu_i^2)
\;=\; \frac{1}{k} \sum\limits_{i=0}^{k-1} \ln(1 + e^{\ln(p_i
\nu_i^2)}) \\
&\geq& \ln\left( 1 + e^{\frac{1}{k} \sum_{i=0}^{k-1} \ln
(p_i \nu_i^2)} \right) \;=\; \ln\left( 1 + \left[
\prod\limits_{i=0}^{k-1} p_i \nu_i^2 \right]^{1/k} \right).
\ea
\eeq
But, for all $0 \leq i \leq k - 1$, we have $\nu_i^2 \geq 
\frac{1}{2} \frac{\la (G_i - A) G_{i+1}^{-1} (G_i - A) u_i,
u_i \ra}{\la G_i u_i, u_i \ra} = \frac{1} {2} \frac{g_
{i+1}^2}{g_i^2}$ by Lemma~\ref{lm-ch-met}, 
\eqref{quad-lin}, and since $G_i u_i = -\nabla f (x_i)$, $A
u_i = \nabla f(x_{i+1}) - \nabla f(x_i)$. Hence, $\prod_
{i=0}^ {k-1} \nu_i^2 \geq \frac{1}{2^k} \frac{g_k^2}
{g_0^2}$, and so $\frac{n}{k} \ln \frac{L}{\mu} \refGE
{sup-prel-quad-2} \ln\left(1 + \frac{1}{2} \left[
\prod_{i=0}^{k-1} p_i \right]^{1/k} \left[ \frac{g_k}
{g_0} \right]^{2/k} \right)$. Rearranging, we obtain $g_k
\leq \left[ \frac{2}{\prod_{i=0}^{k-1} p_i^{1/k}} (
e^{\frac{n}{k} \ln \frac{L} {\mu}} - 1 ) \right]^{k/2} g_0$.
It remains to note that $\lambda_k \leq \sqrt{\frac{L}
{\mu}} \cdot g_k$ and $g_0 \leq \lambda_0$ in view of
\eqref{quad-lin}.\qed
\end{proof}

\BR
As can be seen from \eqref{sup-prel-quad}, the factor $n
\ln \frac{L}{\mu}$ in \eqref{sup-quad} can be improved up to
$\ln\Det(A^{-1}, L B) = \sum_{i=1}^n \ln \frac{L}
{\lambda_i}$, where $\lambda_1, \dots, \lambda_n$ are the
eigenvalues of $A$ relative to $B$. This improved factor can
be significantly smaller than the original one if the
majority of the eigenvalues $\lambda_i$ are much larger than
$\mu$.
\ER

Let us briefly present another approach, which is based on
the \emph{augmented} log-det barrier. The resulting
efficiency estimate will be the same as in
Theorem~\ref{th-sup-quad} up to a slightly worse absolute
constant under the exponent. However, this proof can be
extended onto general nonlinear functions.

\BT\label{th-sup-quad-psi}
In scheme \eqref{qn-quad}, for all $k \geq 1$, we have
$$
\ba{rcl}
\lambda_k &\leq& \left[ \frac{2}{\prod_{i=0}^{k-1} (\tau_i
\frac{\mu}{L} + 1 - \tau_i)^{1/k}} \left (e^{
\frac{13}{6} \frac{n}{k} \ln \frac{L}{\mu}} - 1 \right)
\right]^{k/2} \sqrt{\frac{L}{\mu}} \cdot \lambda_0.
\ea
$$
\ET

\begin{proof}
Without loss of generality, we can assume that $u_i \neq 0$
for all $0 \leq i \leq k$. Denote $\psi_i \Def \psi(G_i,
A)$, $\nu_i \Def \nu(A, G_i, u_i)$, $p_i = \tau_i
\frac{\mu}{L} + 1 - \tau_i$, $g_i
\Def \| \nabla f(x_i) \|_{G_i}^*$ for all $0 \leq i \leq k$.
By Lemma~\ref{lm-prog-psi} and \eqref{quad-lin}, for all $0
\leq i \leq k - 1$, we have $\frac{6}{13} \ln(1 + p_i
\nu_i^2) \leq \psi_i - \psi_{i+1}$. Hence,
\beq\label{sup-2-prel}
\ba{rcl}
\frac{6}{13} \sum\limits_{i=0}^{k-1} \ln(1 + p_i \nu_i^2)
&\leq& \psi_0 - \psi_k \;\refLE{psi-breg}\; \psi_0
\;\refEQ
{qn-quad}\; \psi(L B, A) \\
&\refEQ{def-psi}& \ln\Det(A^{-1}, L B) - \la \frac{1}{L} B^
{-1}, L B - A \ra \;\refLE{quad-mu-L}\; n \ln\frac{L}{\mu},
\ea
\eeq
and we can continue exactly as in the proof of
Theorem~\ref{th-sup-quad}.\qed
\end{proof}

\section{Minimization of General Functions}\label{sec-gen}

In this section, we consider the general unconstrained
minimization problem:
\beq\label{prob-gen}
\ba{rcl}
\min\limits_{x \in \E} f(x),
\ea
\eeq
where $f : \E \to \R$ is a twice continuously differentiable
function with positive definite second derivative. Our goal
is to study the convergence properties of the following
standard quasi-Newton scheme for solving \eqref{prob-gen}:
\beq\label{sch-qn}
\bal \hline
\textbf{Initialization:} Choose $x_0 \in \E$. Set $G_0 = L
B$. \\ \hline
\textbf{For $k \geq 0$ iterate:} \\
1. Update $x_{k+1} = x_k - G_k^{-1} \nabla f(x_k)$. \\
2. Set $u_k = x_{k+1} - x_k$ and choose $\tau_k \in [0, 1]$.
\\
3. Denote $J_k = \int_0^1 \nabla^2 f(x_k + t u_k) dt$.\\
4. Set $G_{k+1} = \Broyd_{\tau_k}(J_k, G_k, u_k)$. \\
\hline
\eal
\eeq
Here $B : \E \to \E^*$ is a self-adjoint positive
definite linear operator, and $L$ is a positive
constant, which together define the initial Hessian
approximation $G_0$.

We assume that there exist constants $\mu > 0$ and $M \geq
0$, such that
\beq\label{mu-L}
\ba{rclrcl}
\mu B &\preceq& \nabla^2 f(x) &\preceq& L B,
\ea
\eeq
\beq\label{sscf}
\ba{rcl}
\nabla^2 f(y) - \nabla^2 f(x) &\preceq& M \| y - x \|_z
\nabla^2 f(w)
\ea
\eeq
for all $x, y, z, w \in \E$. The first assumption
\eqref{mu-L} specifies that, relative to the operator~$B$,
the objective function $f$ is $\mu$-\emph{strongly convex}
and its gradient is $L$-\emph{Lipschitz continuous}. The
second assumption \eqref{sscf} means that $f$ is
$M$-\emph{strongly self-concordant}. This assumption was
recently introduced in \cite{RodomanovNesterov2020a} as a
convenient affine-invariant alternative to the standard
assumption of the Lipschitz second derivative, and is
satisfied at least for any strongly convex function with
Lipschitz continuous Hessian (see \cite[Example~4.1]
{RodomanovNesterov2020a}). The main facts, which we use
about strongly self-concordant functions, are summarized in
the following lemma (see \cite[Lemma~4.1]
{RodomanovNesterov2020a}):
\BL\label{lm-hess-scf}
For any $x, y \in \E$, $J \Def \int_0^1 \nabla^2 f(x + t (y
- x)) dt$, $r \Def \| y - x \|_x$:
\beq\label{hess-J-x}
\ba{rclrcl}
\left( 1 + \frac{M r}{2} \right)^{-1} \nabla^2 f(x)
&\preceq& J &\preceq& \left( 1 + \frac{M r}{2} \right)
\nabla^2 f(x),
\ea
\eeq
\beq\label{hess-J-y}
\ba{rclrcl}
\left( 1 + \frac{M r}{2} \right)^{-1} \nabla^2 f(y)
&\preceq& J &\preceq& \left( 1 + \frac{M r}{2} \right)
\nabla^2 f(y).
\ea
\eeq
\EL
Note that for a quadratic function, we have $M = 0$.

For measuring the convergence rate of
\eqref{sch-qn}, we use the local gradient norm:
\beq\label{def-lam}
\ba{rcl}
\lambda_k &\Def& \| \nabla f(x_k) \|_{x_k}^* \;\refEQ
{def-norms}\; \la \nabla f(x_k), \nabla^2 f(x_k)^{-1} \nabla
f(x_k) \ra^{1/2}.
\ea
\eeq

The local convergence analysis of the scheme~\eqref{sch-qn}
is, in general, the same as the corresponding analysis in
the quadratic case. However, it is much more technical due
to the fact that, in the nonlinear case, the Hessian is no
longer constant. This causes a few problems.

First, there are now several different ways how one can
treat the Hessian approximation $G_k$. One can view it as an
approximation to the Hessian $\nabla^2 f(x_k)$ at the
current iterate $x_k$, to the Hessian $\nabla^2 f(x^*)$ at
the minimizer $x^*$, to the integral Hessian $J_k$ etc. Of
course, locally, due to strong self-concordancy, all these
variants are equivalent since the corresponding Hessians are
close to each other. Nevertheless, from the viewpoint of
technical simplicity of the analysis, some options are
slightly more preferable than others. We find it to be the
most convenient to always think of $G_k$ as an approximation
to the integral Hessian $J_k$.

The second issue is as follows. Suppose we already know what
is the connection between our current Hessian approximation
$G_k$ and the actual integral Hessian $J_k$, e.g., in terms
of the relative eigenvalues and the value of the augmented
log-det barrier potential function \eqref{def-psi}.
Naturally, we want to know how these quantities change after
we update $G_k$ into $G_{k+1}$ at Step~4 of the
scheme~\eqref{sch-qn}. For this, we apply
Lemma~\ref{lm-eigs} and Lemma~\ref{lm-prog-psi}
respectively. However, the problem is that both of these
lemmas will provide us only with the information on the
connection between the update result $G_{k+1}$ and the
\emph{current} integral Hessian $J_k$ (which was used for
performing the update), not the next one $J_{k+1}$.
Therefore, we need to additionally take into account the
errors, resulting from approximating $J_{k+1}$ by $J_k$.

For estimating the errors, which accumulate as a result of
approximating one Hessian by another, it is convenient to
introduce the following quantities\footnote{We follow the
standard convention that the sum over the empty set is
defined as 0, so $\xi_0 = 1$. Similarly, the product over
the empty set is defined as 1.}:
\beq\label{def-r-xi}
\ba{rclrcl}
r_k &\Def& \| u_k \|_{x_k},
\qquad
\xi_k &\Def& e^{M \sum_{i=0}^{k-1} r_i} \quad (\;\geq\; 1),
\qquad k \geq 0.
\ea
\eeq

\BR
The general framework of our analysis is the same as in the
previous paper \cite{RodomanovNesterov2020b}. The main
difference is that now another potential function is used
for establishing the rate of superlinear convergence (Lemma
5.4). However, in order to properly incorporate the new
potential function into the analysis, many parts in the
proof had to be appropriately modified, most notably the
part, related to estimating the region of local convergence.
In any case, the analysis, presented below, is fully
self-contained, and does not require the reader first go
through \cite{RodomanovNesterov2020b}.
\ER

We analyze the method \eqref{sch-qn} in several steps.
The first step is to establish the bounds on the relative
eigenvalues of the Hessian approximations with respect to
the corresponding Hessians.

\BL\label{lm-op-xi}
For all $k \geq 0$, we have
\beq\label{op-hess-xi}
\ba{rclrcl}
\frac{1}{\xi_k} \nabla^2 f(x_k) &\preceq& G_k &\preceq&
\xi_k \frac{L}{\mu} \nabla^2 f(x_k),
\ea
\eeq
\beq\label{op-int-xi}
\ba{rclrcl}
\frac{1}{\xi_{k+1}} J_k &\preceq& G_k &\preceq& \xi_{k+1}
\frac{L}{\mu} J_k.
\ea
\eeq
\EL

\begin{proof}
For $k=0$, \eqref{op-hess-xi} follows from \eqref{mu-L} and
the fact that $G_0 = L B$ and $\xi_0 = 1$. Now suppose that
$k \geq 0$, and that \eqref{op-hess-xi} has
already been proved for all indices up to~$k$. Then,
applying Lemma~\ref{lm-hess-scf} to \eqref{op-hess-xi}, we
obtain
\beq\label{op-int-aux}
\ba{rclrcl}
\frac{1}{\xi_k \left( 1 + \frac{M r_k}{2} \right)} J_k
&\preceq& G_k &\preceq& \left(1 + \frac{M r_k}{2} \right)
\xi_k \frac{L}{\mu} J_k.
\ea
\eeq
Since $(1 + \frac{M r_k}{2}) \xi_k \leq \xi_{k+1}$ by 
\eqref{def-r-xi}, this proves \eqref{op-int-xi} for the
index $k$. Applying Lemma~\ref{lm-eigs} to
\eqref{op-int-aux}, we get $\frac{1}{\xi_k( 1 + \frac{M
r_k}{2})} J_k \preceq G_{k+1}
\preceq (1 + \frac{M r_k}{2} ) \xi_k \frac{L}{\mu} J_k$, and
so
$$
\ba{rcl}
G_{k+1} &\refPE{hess-J-y}& \left( 1 + \frac{M r_k}{2}
\right)^2 \xi_k \frac{L}{\mu} \nabla^2 f(x_{k+1}) \;\refPE
{def-r-xi}\; \xi_{k+1} \frac{L}{\mu} \nabla^2 f(x_{k+1}), \\
G_{k+1} &\refSE{hess-J-y}& \frac{1}{\left( 1 + \frac{M r_k}
{2} \right)^2 \xi_k} \nabla^2 f(x_{k+1}) \;\refSE{def-r-xi}\;
\frac{1}{\xi_{k+1}} \nabla^2 f(x_{k+1}).
\ea
$$
This proves \eqref{op-hess-xi} for the index $k+1$, and we
can continue by induction.\qed
\end{proof}

\BC
For all $k \geq 0$, we have
\beq\label{r-ubd}
\ba{rcl}
r_k &\leq& \xi_k \lambda_k.
\ea
\eeq
\EC

\begin{proof}
Indeed,
\[
\ba[b]{rcl}
r_k \;\refEQ{def-r-xi}\; \| u_k \|_{x_k} &\refEQ{sch-qn}&
\la \nabla f(x_k), G_k^{-1} \nabla^2 f(x_k) G_k^{-1} \nabla
f(x_k) \ra^{1/2} \\
&\refLE{op-hess-xi}& \xi_k \la \nabla f(x_k), \nabla^2 f 
(x_k)^{-1} \nabla f(x_k) \ra^{1/2} \;\refEQ{def-lam}\; \xi_k
\lambda_k.
\ea\qedEA
\]
\end{proof}

The second step in our analysis is to establish a
preliminary version of the linear convergence theorem for
the scheme \eqref{sch-qn}.

\BL\label{lm-lin-xi}
For all $k \geq 0$, we have
\beq\label{lam-lin-xi}
\ba{rcl}
\lambda_k &\leq& \sqrt{\xi_k} \lambda_0 \prod\limits_{i=0}^
{k-1} q_i,
\ea
\eeq
where
\beq\label{def-q}
\ba{rcl}
q_i &\Def& \max\left\{ 1 - \frac{\mu}{\xi_{i+1} L}, \xi_
{i+1} - 1 \right\}.
\ea
\eeq
\EL

\begin{proof}
Let $k, i \geq 0$ be arbitrary. By Taylor's formula, we have
\beq\label{lin-aux}
\ba{rcl}
\nabla f(x_{i+1}) &\refEQ{sch-qn}& \nabla f(x_i) + J_i u_i
\;\refEQ{sch-qn}\; J_i (J_i^{-1} - G_i^{-1}) \nabla f(x_i).
\ea
\eeq
Hence,
\beq\label{lin-lam-prel}
\ba{rcl}
\lambda_{i+1} &\refEQ{def-lam}& \la \nabla f(x_{i+1}),
\nabla^2 f(x_{i+1})^{-1} \nabla f(x_{i+1}) \ra^{1/2} \\
&\refLE{hess-J-y}& \sqrt{1 + \frac{M r_i}{2}} \la \nabla f
(x_{i+1}), J_i^{-1} \nabla f(x_{i+1}) \ra^{1/2} \\
&\refEQ{lin-aux}& \sqrt{1 + \frac{M r_i}{2}} \la \nabla f
(x_i), (J_i^{-1} - G_i^{-1}) J_i (J_i^{-1} - G_i^{-1})
\nabla f(x_i) \ra^{1/2}.
\ea
\eeq
Note that $-(\xi_{i+1} - 1) J_i^{-1} \refPE{op-int-xi}
J_i^{-1} - G_i^{-1} \refPE{op-int-xi} \left( 1 - \frac{\mu}
{\xi_{i+1} L} \right) J_i^{-1}$. Therefore,
$$
\ba{rcl}
(J_i^{-1} - G_i^{-1}) J_i (J_i^{-1} - G_i^{-1}) &\refPE
{def-q}& q_i^2 J_i^{-1} \;\refPE{hess-J-x}\; q_i^2 \left( 1
+ \frac{M r_i}{2} \right) \nabla^2 f(x_i)^{-1}.
\ea
$$
Thus, $\lambda_{i+1} \leq \left( 1 + \frac{M r_i}{2}\right)
q_i \lambda_i$ in view of \eqref{lin-lam-prel} and 
\eqref{def-lam}. Consequently,
\[
\ba[b]{rcl}
\lambda_k &\leq& \lambda_0 \prod\limits_{i=0}^{k-1}
\left( 1 + \frac{M r_i}{2} \right) q_i \;\leq\; \lambda_0
\prod\limits_{i=0}^{k-1} e^{\frac{M r_i}{2}} q_i \;\refEQ
{def-r-xi}\; \sqrt{\xi_k} \lambda_0 \prod\limits_{i=0}^{k-1}
q_i.
\ea\qedEA
\]
\end{proof}

Next, we establish a preliminary version of the theorem on
superlinear convergence of the scheme \eqref{sch-qn}.
The proof uses the augmented log-det barrier potential
function and is essentially a generalization of the
corresponding proof of Theorem~\ref{th-sup-quad-psi}.

\BL\label{lm-sup-xi}
For all $k \geq 1$, we have
\beq\label{lam-sup-xi}
\ba{rcl}
\lambda_k &\leq& \left[ \frac{1 + \xi_k}
{\prod_{i=0}^{k-1} (\tau_i \frac{\mu}{\xi_{i+1}^2 L} + 1 -
\tau_i)^{1/k}} \left( e^{\frac{13}{6} \frac{n}{k} \ln \left(
\xi_{k+1}^{\xi_{k+1}} \frac{L}{\mu} \right)} - 1 \right)
\right]^{k/2} \sqrt{\xi_k \frac{L}{\mu}}
\cdot \lambda_0.
\ea
\eeq
\EL

\begin{proof}
Without loss of generality, assume that $u_i \neq 0$
for all $0 \leq i \leq k$. Denote $\psi_i \Def \psi(G_i,
J_i)$, $\tilde{\psi}_{i+1}
\Def \psi(G_{i+1}, J_i)$, $\nu_i \Def \nu(J_i, G_i, u_i)$,
$p_i \Def \tau_i \frac{\mu}{\xi_{i+1}^2 L} + 1 - \tau_i$,
and $g_i \Def \| \nabla f(x_i) \|_{G_i}^*$ for any $0 \leq i
\leq k$.

Let $0 \leq i \leq k - 1$ be arbitrary. By
Lemma~\ref{lm-prog-psi} and \eqref{op-int-xi}, we have
\beq\label{psi-prog-prel}
\ba{rcl}
\frac{6}{13} \ln\left( 1 + p_i \nu_i^2 \right) &\leq& \psi_i
- \tilde{\psi}_{i+1} \;=\; \psi_i - \psi_{i+1} + \Delta_i,
\ea
\eeq
where
\beq\label{def-Delta}
\ba{rcl}
\Delta_i &\Def& \psi_{i+1} - \tilde{\psi}_{i+1} \;\refEQ
{def-psi}\; \la G_{i+1}^{-1}, J_{i+1} - J_i \ra + \ln\Det
(J_{i+1}^{-1}, J_i).
\ea
\eeq
Note that $J_i \succeq (1 + \frac{M r_i}{2})^{-1}
\nabla^2 f(x_{i+1}) \succeq (1 + \frac{M r_i}{2})^
{-1} (1 + \frac{M r_{i+1}}{2})^{-1} J_{i+1}$ in view of
\eqref{hess-J-y} and \eqref{hess-J-x}. In
particular,
$J_i \succeq e^{-\frac{M}{2} (r_i + r_{i+1})} J_{i+1}
\succeq (1 - \frac{M}{2} (r_i + r_{i+1})) J_{i+1}$.
Therefore, $J_{i+1} - J_i \preceq \frac{M}{2} (r_i + r_
{i+1}) J_{i+1}$, and
so
$$
\ba{rcl}
\sum\limits_{i=0}^{k-1} \la G_{i+1}^{-1}, J_{i+1} - J_i \ra
&\leq& \frac{M}{2} \sum\limits_{i=0}^{k-1} (r_i + r_{i+1})
\la G_{i+1}^{-1}, J_{i+1} \ra
\\
&\refLE{op-int-xi}& n \frac{M}{2} \sum\limits_{i=0}^{k-1}
\xi_{i+2} (r_i + r_{i+1}) \;\refLE{def-r-xi}\; n \xi_{k+1}
\frac{M}{2} \sum\limits_{i=0}^{k-1} (r_i + r_{i+1}) \\
&\leq& n \xi_{k+1} M \sum\limits_{i=0}^k
r_i \;\refEQ {def-r-xi}\; n \xi_{k+1} \ln \xi_{k+1}.
\ea
$$
Consequently,
\beq\label{sum-Delta}
\ba{rcl}
\sum\limits_{i=0}^{k-1} \Delta_i &\refLE{def-Delta}& n \xi_
{k+1} \ln \xi_{k+1} + \ln\Det(J_k^{-1}, J_0).
\ea
\eeq
Summing up \eqref{psi-prog-prel}, we thus obtain
$$
\ba{rcl}
&&\hspace{-5em}\frac{6}{13} \sum\limits_{i=0}^{k-1} \ln(1 +
p_i \nu_i^2)
\leq \psi_0 - \psi_k + \sum\limits_{i=0}^{k-1} \Delta_i
\;\refLE{psi-breg}\; \psi_0 + \sum\limits_{i=0}^{k-1}
\Delta_i \\
&\refEQ{def-psi}& \ln\Det(J_0^{-1}, L B) - \la \frac{1}{L}
B^{-1}, L B - J_0 \ra + \sum\limits_{i=0}^{k-1} \Delta_i \\
&\refLE{sum-Delta}& \ln\Det(J_k^{-1}, L B) - \la \frac{1}{L}
B^{-1}, L B - J_0 \ra + n \xi_{k+1} \ln \xi_{k+1} \\
&\refLE{mu-L}& n \ln \frac{L}{\mu} + n \xi_{k+1}
\ln \xi_{k+1} \;=\; n \ln \left( \xi_{k+1}^{\xi_{k+1}} 
\frac{L}{\mu} \right).
\ea
$$
By the convexity of function $t \mapsto \ln(1 +
e^t)$, it follows that
\beq\label{sup-prel}
\ba{rcl}
&&\hspace{-3em}\frac{13}{6} \frac{n}{k} \ln \left( \xi_
{k+1}^{\xi_{k+1}} \frac{L}{\mu} \right) \;\geq\; \frac{1}{k}
\sum\limits_{i=0}^{k-1} \ln(1 +
p_i \nu_i^2) \;=\; \frac{1}{k} \sum\limits_{i=0}^{k-1} \ln(1
+ e^{\ln(p_i \nu_i^2)}) \\
&\geq& \ln\left( 1 + e^{\frac{1}{k} \sum_{i=0}^ {k-1} \ln
(p_i \nu_i^2)} \right) \;=\; \ln\left( 1 + \left[
\prod\limits_{i=0}^{k-1} p_i \nu_i^2 \right]^{1/k} \right).
\ea
\eeq
At the same time, $\nu_i^2 \geq \frac{1}{1 +
\xi_{i+1}} \frac{\la (G_i - J_i) G_{i+1}^{-1} (G_i - J_i)
u_i, u_i \ra}{\la G_i u_i, u_i \ra} = \frac{1}{1 + \xi_
{i+1}} \frac{g_{i+1}^2}{g_i^2}$ in view of
Lemma~\ref{lm-ch-met}, \eqref{op-int-xi} and since $G_i u_i
= -\nabla f(x_i)$, $J_i u_i = \nabla f(x_ {i+1}) - \nabla
f(x_i)$. Hence, we can write $\prod_{i=0}^{k-1} \nu_i^2 \geq
\frac{g_k^2}{g_0^2} \prod_{i=0}^{k-1} \frac{1}{1 +
\xi_{i+1}} \refGE{def-r-xi} \frac{1}{(1 + \xi_k)^k} 
\frac{g_k^2}{g_0^2}$. Consequently, $\frac{13}{6} 
\frac{n}{k} \ln( \xi_{k+1}^{\xi_{k+1}} \frac{L}{\mu})
\refGE{sup-prel} \ln\left( 1 + \frac{\prod_{i=0}^{k-1} p_i^
{1/k}}{1 + \xi_k} \left[ \frac{g_k}{g_0} \right]^{2/k}
\right)$. Rearranging, we obtain that $g_k \leq \left[ 
\frac{1 + \xi_k}{\prod_{i=0}^{k-1} p_i^{1/k}} (e^{\frac{13}
{6} \frac{n}{k} \ln ( \xi_ {k+1}^{\xi_{k+1}} \frac{L}{\mu}
)} - 1 ) \right]^{k/2} g_0$. But $\lambda_k \leq \sqrt{\xi_k
\frac{L}{\mu}} \cdot g_k$ by \eqref{op-hess-xi}, and
$g_0 \leq \lambda_0$ in view of \eqref{mu-L} and the fact
that $G_0 = L B$.\qed
\end{proof}

In the quadratic case ($M = 0$), we have $\xi_k
\equiv 1$ (see \eqref{def-r-xi}), and Lemmas~\ref{lm-op-xi}
and \ref{lm-lin-xi} reduce to the
already known Theorem~\ref{th-lin-quad}, and
Lemma~\ref{lm-sup-xi} reduces to the already known
Theorem~\ref{th-sup-quad}. In the general case, the
quantities $\xi_k$ can grow with iterations. However, as
we will see in a moment, by requiring the initial
point~$x_0$ in the scheme \eqref{sch-qn} to be sufficiently
close to the solution, we can still ensure that $\xi_k$
stay \emph{uniformly bounded} by a sufficiently small
absolute constant. This allows us to recover all the main
results of the quadratic case.

To write down the region of local convergence
of~\eqref{sch-qn}, we need to introduce one more quantity,
related to the starting moment of superlinear
convergence\footnote{Hereinafter, $\lceil t \rceil$ for $t
> 0$ denotes the smallest positive integer greater or equal
to $t$.}:
\beq\label{def-K0-tau}
\ba{rclrcl}
K_0 &\Def& \left\lceil \frac{1}{\tau \frac{4 \mu}{9 L} + 1 -
\tau} 8 n \ln \frac{2 L}{\mu} \right\rceil,
\qquad
\tau &\Def& \sup\limits_{k \geq 0} \tau_k \quad (\;\leq\;
1).
\ea
\eeq
For DFP ($\tau_k \equiv 1$) and BFGS ($\tau_k \equiv 0)$, we
have respectively
\beq\label{K0-DFP-BFGS}
\ba{rclrcl}
K_0^{\DFP} &=& \left\lceil \frac{18 n L}{\mu} \ln \frac{2 L}
{\mu} \right\rceil,
\qquad
K_0^{\BFGS} &=& \left\lceil 8 n \ln \frac{2 L}{\mu}
\right\rceil.
\ea
\eeq

Now we are ready to prove the main result of this section.

\BT\label{th-sup}
Suppose that, in scheme \eqref{sch-qn}, we have
\beq\label{lam-ini}
\ba{rcl}
M \lambda_0 &\leq& \frac{\ln \frac{3}{2}}{\left(
\frac{3}{2} \right)^{\frac{3}{2}}} \max\left\{
\frac{\mu}{2 L}, \frac{1}{K_0 + 9} \right\}.
\ea
\eeq
Then, for all $k \geq 0$,
\beq\label{op-hess}
\ba{rclrcl}
\frac{2}{3} \nabla^2 f(x_k) &\preceq& G_k &\preceq&
\frac{3 L}{2 \mu} \nabla^2 f(x_k),
\ea
\eeq
\beq\label{lam-lin}
\ba{rclrcl}
\lambda_k &\leq& \left( 1 - \frac{\mu}{2 L} \right)^k
\sqrt{\frac{3}{2}} \cdot \lambda_0,
\ea
\eeq
and, for all $k \geq 1$,
\beq\label{lam-sup}
\ba{rcl}
\lambda_k &\leq& \left[\frac{5}{2 \prod_{i=0}^{k-1} (\tau_i
\frac{4 \mu}{9 L} + 1 - \tau_i)^{1/k}} \left( e^{\frac{13}
{6} \frac{n}{k} \ln \frac{2 L}{\mu} } - 1
\right) \right]^{k/2} \sqrt{\frac{3 L}{2 \mu}} \cdot
\lambda_0.
\ea
\eeq
\ET

\begin{proof}
Let us prove by induction that, for all $k \geq 0$, we have
\beq\label{xi-ubd}
\ba{rcl}
\xi_k &\leq& \frac{3}{2}.
\ea
\eeq
Clearly, \eqref{xi-ubd} is satisfied for $k=0$
since $\xi_0 = 1$. It is also satisfied for $k=1$
since $\xi_1 \refEQ{def-r-xi} e^{M r_0} \refLE{r-ubd}
e^{\xi_0 M \lambda_0} \refEQ{def-r-xi} e^{M \lambda_0}
\refLE{lam-ini} \frac{3}{2}$.

Now let $k \geq 0$, and suppose that \eqref{xi-ubd} has
already been proved for all indices up to~$k+1$. Then,
applying Lemma~\ref{lm-op-xi}, we obtain \eqref{op-hess} for
all indices up to $k+1$. Applying now Lemma~\ref{lm-lin-xi}
and using for all $0 \leq i \leq k$ the relation $q_i
\refEQ{def-q} \max\{ 1 - \frac{\mu}{\xi_{i+1} L},
\xi_{i+1} - 1 \} \refLE{xi-ubd} \max\{1 - \frac{2 \mu}{3 L},
\frac{1}{2} \} \leq 1 - \frac{\mu}{2 L}$, we obtain 
\eqref{lam-lin} for all indices up to $k+1$.
Finally, if $k \geq 1$, then, applying Lemma~\ref{lm-sup-xi}
and using that $\xi_{i+1}^{\xi_{i+1}} \refLE{xi-ubd} ( 
\frac{3}{2} )^{\frac{3}{2}} = \frac{3}{2} \sqrt{\frac{3}
{2}} \leq \frac{3}{2} (1 + \frac{1}{4}) = \frac{15}{8} \leq
2$ for all $0 \leq i \leq k$, we obtain \eqref{lam-sup} for
all indices up to $k$. Thus, at this moment, \eqref{op-hess}
and~\eqref{lam-lin} are proved for all indices up to $k+1$,
while \eqref{lam-sup} is proved only up to $k$.

To finish the inductive step, it remains to prove that
\eqref{xi-ubd} is satisfied for the index $k+2$,
or, equivalently, in view of \eqref{def-r-xi}, that $M
\sum_{i=0}^{k+1} r_i \leq \ln \frac{3}{2}$. Since
$M \sum_{i=0}^{k+1} r_i \leq M
\sum_{i=0}^{k+1} \xi_i \lambda_i \leq 
\frac{3}{2} M \sum_{i=0}^{k+1} \lambda_i$ in view of 
\eqref{r-ubd}
and \eqref{xi-ubd} respectively, it suffices to show that
$\frac{3}{2} M \sum_{i=0}^{k+1} \lambda_i \leq \ln
\frac{3}{2}$.

Note that
\beq\label{sum-lam-1}
\ba{rcl}
\frac{3}{2} M \sum\limits_{i=0}^{k+1} \lambda_i
&\refLE{lam-lin}& \left( \frac{3}{2} \right)^{\frac{3}{2}} M
\lambda_0 \sum\limits_{i=0}^{k+1} \left( 1 - \frac{\mu}{2 L}
\right)^i \;\leq\; \left( \frac{3}{2} \right)^{\frac{3}{2}} 
\frac{2
L} {\mu} M \lambda_0.
\ea
\eeq
Therefore, if we could prove that
\beq\label{sum-lam-2}
\ba{rcl}
\frac{3}{2} M \sum\limits_{i=0}^{k+1} \lambda_i
&\leq& \left( \frac{3}{2} \right)^{\frac{3}{2}} (K_0 + 9) M
\lambda_0,
\ea
\eeq
then, combining \eqref{sum-lam-1} and \eqref{sum-lam-2}, we
would obtain
$$
\ba{rcl}
\frac{3}{2} M \sum\limits_{i=0}^{k+1} \lambda_i &\leq& \left
(\frac{3}{2} \right)^{\frac{3}{2}} \min\left \{ \frac{2 L}
{\mu}, K_0 + 9 \right\} M \lambda_0 \;\refLE{lam-ini}\; \ln 
\frac{3}
{2},
\ea
$$
which is exactly what we need. Let us prove 
\eqref{sum-lam-2}. If $k \leq K_0$, in view of 
\eqref{lam-lin}, we have $\frac{3} {2} M
\sum_{i=0}^{k+1} \lambda_i \leq \left(
\frac{3}{2} \right)^{\frac{3}{2}} (k + 2) M \lambda_0
\leq \left( \frac{3}{2} \right)^{\frac{3} {2}} (K_0 + 2)
M \lambda_0$, and \eqref{sum-lam-2} follows. Therefore, from
now on, we
can assume that $k \geq K_0$. Then\footnote{We will estimate
the second sum using \eqref{lam-sup}. However, recall that,
at this moment, \eqref{lam-sup} is proved only up to the
index $k$. This is the reason why we move $\lambda_{k+1}$
into the first sum.},
$$
\ba{rcl}
\frac{3}{2} M \sum\limits_{i=0}^{k+1} \lambda_i &=& \frac{3}
{2} M \left( \sum\limits_{i=0}^{K_0-1} \lambda_i + \lambda_
{k+1} \right) + \frac{3}{2} M \sum\limits_{i=K_0}^k
\lambda_i \\
&\refLE{lam-lin}& \left( \frac{3}{2} \right)^{\frac{3}{2}} 
(K_0 + 1) M \lambda_0 + \frac{3}{2} M \sum\limits_{i=K_0}^k
\lambda_i.
\ea
$$
It remains to show $\frac{3}{2} M \sum_{i=K_0}^k \lambda_i
\leq \left( \frac{3}{2} \right)^{\frac{3}{2}} 8 M
\lambda_0$. We can do this using \eqref{lam-sup}.

First, let us make some estimations. Clearly, for all $0 < t
< 1$, we have $e^t = \sum_{j=0}^\infty \frac{t^j} {j!} \leq
1 + t + \frac{t^2}{2} \sum_{j=0}^\infty t^j = 1 + t ( 1
+ \frac{t}{2 (1 - t)})$. Hence, for
all $0 < t \leq 1$, we obtain $e^{\frac{13 t}{48}}
- 1 \leq \frac{13 t}{48} ( 1 + \frac{\frac{13}{48}}{2 
(1 - \frac{13}{48})} ) = \frac{13 t}{48} \cdot 
\frac{83}{70} \leq \frac{13 t}{48} \cdot \frac{6}{5} = 
\frac{13 t}{40}$, and so
\beq\label{exp-aux-2}
\ba{rcl}
\left[ \frac{5}{2 t} \left( e^{\frac{13 t}{48}} - 1
\right) \right]^{1/2} &\leq& \sqrt{ \frac{5}
{2 t} \cdot \frac{13 t} {40}} \;=\; \sqrt{\frac{13}{16}}
\;\leq\; \frac{11}{12}.
\ea
\eeq
At the same time, $\frac{11}{12} = 1 - \frac{1}{12} \leq
e^{-\frac{1}{12}}$. Hence,
\beq\label{sup-aux-2}
\ba{rcl}
\left( \frac{11}{12} \right)^{K_0} \sqrt{\frac{L}{\mu}}
&\refLE{def-K0-tau}& \left( \frac{11}{12} \right)^{8 
\ln \frac{2 L}{\mu}} \sqrt{\frac{L}{\mu}} \;\leq\; e^
{-\frac{2}{3} \ln \frac{2 L}{\mu}} \sqrt{\frac{L}{\mu}}
\;=\; \left( \frac{2 L}{\mu} \right)^{-\frac{2}{3}} \sqrt{
\frac{L}{\mu}} \\
&=& 2^{-\frac{2}{3}} \left( \frac{L}{\mu} \right)^{-
\frac{1}{6}} \;\leq\; 2^{-\frac{2}{3}} \;\leq\;
\frac{2}{3}.\vspace{-1em}
\ea
\eeq
Thus, for all $K_0 \leq i \leq k$, and $p \Def
\tau
\frac{4 \mu}{9 L} + 1 - \tau \;\refLE {def-K0-tau}\; \prod_
{j=0}^{i-1} (\tau_i \frac{4 \mu}{9 L} + 1 - \tau_i)^{1/i}$:
$$
\ba{rcl}
\lambda_i &\refLE{lam-sup}& \left[ \frac{5}{2
p} \left( e^{ \frac{13}{6} \frac{n}{i} \ln \frac{2 L}{\mu} }
- 1 \right) \right]^{i/2} \sqrt{\frac{3 L}{2 \mu}} \cdot 
\lambda_0 \\
&\refLE{def-K0-tau}& \left[ \frac{5}{2 p} \left( e^{\frac{13
p}{48}} - 1 \right) \right]^{i/2} \sqrt{\frac{3 L}{2 \mu}}
\cdot \lambda_0 \;\refLE{exp-aux-2}\; \left( \frac{11}{12}
\right)^i \sqrt{\frac{3 L}{2 \mu}} \cdot \lambda_0 \\
&=& \left( \frac{11}{12} \right)^{i-K_0} \left( \frac{11}
{12} \right)^{K_0} \sqrt{\frac{3 L}{2 \mu}} \cdot \lambda_0
\;\refLE{sup-aux-2}\; \left( \frac{11}{12} \right)^{i-K_0}
\frac{2}{3} \cdot \sqrt{\frac{3}{2}} \cdot \lambda_0.
\ea
$$
Hence, $\frac{3}{2} M \sum_{i=K_0}^k
\lambda_i \leq (\frac{3}{2})^{\frac{3}{2}} M
\lambda_0 \cdot \frac{2}{3} \sum_{i=K_0}^k (\frac{11}{12})^
{i-K_0} \leq (\frac{3}{2})^{ \frac{3}{2}} 8 M
\lambda_0$.\qed
\end{proof}

\BR
In accordance with Theorem~\ref{th-sup}, the parameter $M$
of strong self-concordancy affects only the size of the
region of local convergence of the process~\eqref{sch-qn},
and not its rate of convergence. We do not know
whether this is an artifact of the analysis or not, but
it might be an interesting topic for future research. For a
quadratic function, we have $M = 0$, and so the
scheme~\eqref{sch-qn} is globally convergent.
\ER

The region of local convergence, specified by
\eqref{lam-ini}, depends on the \emph{maximum} of two
quantities: $\frac{\mu}{L}$ and $\frac{1}{K_0}$. For DFP,
the $\frac{1}{K_0}$ part in this maximum is in fact
redundant, and its region of local convergence is simply
inversely proportional to the condition number: $O\left( 
\frac{\mu}{L} \right)$. However, for BFGS, the $\frac{1}
{K_0}$ part does not disappear, and we obtain the following
region of local convergence:
$$
\ba{rcl}
M \lambda_0 &\leq&
\max\left\{ O\left(\frac{\mu}{L}\right), \ O\left(\frac{1}{n
\ln \frac{2 L}{\mu}} \right) \right\}.
\ea
$$
Clearly, the latter region can be much bigger than the
former when the condition number $\frac{L}{\mu}$ is
significantly larger than the dimension $n$.

\BR
The previous estimate of the size of the region of local
convergence, established in \cite{RodomanovNesterov2020b},
was $O(\frac{\mu}{L})$ for both DFP and BFGS.
\ER

\BE
Consider the functions
\[
\ba{rclrcl}
f(x) &\Def& f_0(x) + \frac{\mu}{2} \| x \|^2,
\qquad
f_0(x) &\Def& \ln\left( \sum\limits_{i=1}^m e^{\la a_i, x \ra
+ b_i} \right), \qquad x \in \E,
\ea
\]
where $a_i \in \E^*$, $b_i \in \R$, $i = 1, \ldots, m$, $\mu
> 0$, and $\| \cdot \|$ is the Euclidean norm, induced by
the operator $B$. Let $\gamma > 0$ be such that
\[
\ba{rcl}
\| a_i \|_* &\leq& \gamma, \qquad i = 1, \ldots, m,
\ea
\]
where $\| \cdot \|_*$ is the norm conjugate to $\| \cdot
\|$. Define
\[
\ba{rcl}
\pi_i(x) &\Def& \frac{e^{\la a_i, x \ra + b_i}}{\sum_{j=1}^m
e^{\la a_j, x \ra + b_j}}, \qquad x \in \E, \quad i = 1,
\ldots, m.
\ea
\]
Clearly, $\sum_{i=1}^m \pi_i(x) = 1$, $\pi_i(x) > 0$ for all
$x \in \E$, $i = 1, \ldots, m$. It is not difficult to check
that, for all $x, h \in \E$, we have\footnote{$D^3
f_0(x)[h, h, h] = \frac{d^3}{d t^3} f_0(x + t h)
\big\vert_{t = 0}$ is the third derivative of $f$ along the
direction $h$.}
\[
\ba{rcl}
\la \nabla f_0(x), h \ra &=& \sum\limits_{i=1}^m \pi_i(x)
\la a_i, h \ra \;\leq\; \gamma.\\
\la \nabla^2 f_0(x) h, h \ra &=& \sum\limits_{i=1}^m
\pi_i(x) \la a_i - \nabla f_0(x), h \ra^2 \\
&=& \sum\limits_{i=1}^m \pi_i(x) \la a_i, h \ra^2 - \la
\nabla f_0(x), h \ra^2 \;\leq\; \gamma^2 \| h \|^2, \\
D^3 f_0(x) [h, h, h] &=& \sum\limits_{i=1}^m \pi_i(x) \la
a_i - \nabla f_0(x), h \ra^3 \\
&\leq& 2 \gamma \| h \| \la \nabla^2 f_0(x) h, h \ra
\;\leq\; 2 \gamma^3 \| h \|^3.
\ea
\]
Thus, $f_0$ is a convex function with $\gamma^2$-Lipschitz
gradient and $(2 \gamma^3)$-Lipschitz Hessian. Consequently,
the function $f$ is $\mu$-strongly convex with $L$-Lipschitz
gradient, $(2 \gamma^3)$-Lipschitz Hessian, and, in view of
\cite[Example~4.1]{RodomanovNesterov2020a}, $M$-strongly
self-concordant, where
\[
\ba{rclrcl}
L &\Def& \gamma^2 + \mu,
\qquad
M &\Def& \frac{2 \gamma^3}{\mu^{3/2}}.
\ea
\]

Let the regularization parameter $\mu$ be sufficiently
small, namely $\mu \leq \gamma^2$. Denote $Q \Def
\frac{\gamma^2}{\mu} \geq 1$. Then, $Q \leq \frac{L}{\mu}
\leq 2 Q$, $M = 2 Q^{3/2}$, so, according to
\eqref{lam-ini}, the region of local convergence of BFGS
can be described as follows:
\[
\ba[b]{rcl}
\lambda_0 &\leq& \max\left\{ O\left( \frac{1}{Q^{5/2}}
\right), O\left( \frac{1}{n Q^{3/2} \ln(4 Q)} \right)
\right\}.
\ea\qedEA
\]
\EE

\section{Discussion}\label{sec-disc}

Let us compare the new convergence rates, obtained in this
paper for the classical DFP and BFGS methods, with the
previously known ones from~\cite{RodomanovNesterov2020b}.
Since the estimates for the general nonlinear case differ
from those for the quadratic one just in absolute constants,
we only discuss the latter case.

In what follows, we use our standard notation: $n$ is the
dimension of the space, $\mu$ is the strong convexity
parameter, $L$ is the Lipschitz constant of the gradient,
and $\lambda_k$ is the local norm of the gradient at
the $k$th iteration.

For BFGS, the previously known rate (see~\cite[Theorem~3.2]
{RodomanovNesterov2020b}) is
\beq\label{bfgs-prev}
\ba{rcl}
\lambda_k &\leq& \left( \frac{n L}{\mu k} \right)^{k/2}
\lambda_0.
\ea
\eeq
Although \eqref{bfgs-prev} is formally valid for all $k \geq 1$,
it becomes useful\footnote{Indeed, according to
Theorem~\ref{th-lin-quad}, we have at least
$\lambda_k \leq (1 - \frac{\mu}{L})^k \lambda_0$ for all
$k \geq 0$.} only after
\beq\label{bfgs-st-prev}
\ba{rcl}
\widehat{K}_0^{\BFGS} &\Def& \frac{n L}{\mu}
\ea
\eeq
iterations. Thus, $\widehat{K}_0^{\BFGS}$ can be thought of
as the \emph{starting moment} of the superlinear
convergence, according to the estimate \eqref{bfgs-prev}.

In this paper, we have obtained a new estimate
(Theorem~\ref{th-sup-quad}):
\beq\label{bfgs-new}
\ba{rcl}
\lambda_k &\leq& \left[ 2 \left( e^{\frac{n}{k} \ln
\frac{L}{\mu}} - 1 \right) \right]^{k/2} \sqrt{\frac{L}
{\mu}} \cdot \lambda_0.
\ea
\eeq
Its starting moment of superlinear convergence can be
described as follows:
\beq\label{bfgs-st-new}
\ba{rcl}
K_0^{\BFGS} &\Def& 4 n \ln \frac{L}{\mu}.
\ea
\eeq
Indeed, since $e^t \leq \frac{1}{1-t} = 1 + \frac{t}{1-t}$
for any $t < 1$, we have, for all $k \geq K_0^{\BFGS}$,
\beq\label{disc-aux-1}
\ba{rcl}
e^{\frac{n}{k} \ln \frac{L}{\mu}} - 1 &\leq&
\frac{\frac{n}{k} \ln \frac{L}{\mu}}{1 - \frac{n}{k} \ln
\frac{L}{\mu}} \;\refLE{bfgs-st-new}\; \frac{\frac{n}{k} \ln
\frac{L}{\mu}}{1 - \frac{1}{4}} \;=\; \frac{4 n}{3 k} \ln
\frac{L}{\mu}.
\ea
\eeq
At the same time, for all $k \geq K_0^{\BFGS}$:
\beq\label{disc-aux-2}
\ba{rcl}
\sqrt{\frac{L}{\mu}} &=& e^{\frac{1}{2} \ln \frac{L}{\mu}}
\;\refLE{bfgs-st-new}\; e^{\frac{k}{8}} \;=\; (e^{\frac{1}
{4}} )^ {k/2} \;\leq\; \left( \frac{4}{3} \right)^{k/2}
\;\leq\; \left( \frac{3}{2} \right)^{k/2}.
\ea
\eeq
Hence, according the new estimate \eqref{bfgs-new}, for all
$k \geq K_0^{\BFGS}$:
\beq\label{bfgs-new-2}
\ba{rcl}
\lambda_k &\refLE{disc-aux-1}& \left( \frac{8 n}{3 k} \ln
\frac{L}{\mu} \right)^{k/2} \sqrt{\frac{L}{\mu}} \cdot
\lambda_0 \;\refLE{disc-aux-2}\; \left( \frac{4 n}{k} \ln
\frac{L}{\mu} \right)^{k/2} \lambda_0 \qquad (\;\refLE
{bfgs-st-new}\; \lambda_0).
\ea
\eeq
Comparing the previously known efficiency estimate
\eqref{bfgs-prev} and its starting moment of superlinear
convergence \eqref{bfgs-st-prev} with the new ones
\eqref{bfgs-new-2}, \eqref{bfgs-st-new}, we thus conclude
that we manage to put the condition number $\frac{L}{\mu}$
\emph{under the logarithm}.

For DFP, the previously known rate (see~\cite
[Theorem~3.2]{RodomanovNesterov2020b}) is
$$
\ba{rcl}
\lambda_k &\leq& \left( \frac{n L^2}{\mu^2 k} \right)^{k/2}
\lambda_0
\ea
$$
with the following starting moment of the superlinear
convergence:
\beq\label{dfp-st-prev}
\ba{rcl}
\widehat{K}_0^{\DFP} &\Def& \frac{n L^2}{\mu^2}.
\ea
\eeq
The new rate, which we have obtained in this paper
(Theorem~\ref{th-sup-quad}), is
\beq\label{dfp-new}
\ba{rcl}
\lambda_k &\leq& \left[ \frac{2 L}{\mu} \left( e^{\frac{n}
{k} \ln \frac{L}{\mu}} - 1 \right) \right]^{k/2} \sqrt{
\frac{L}{\mu}} \cdot \lambda_0.
\ea
\eeq
Repeating the same reasoning as above, we can easily
obtain that the new starting moment of the superlinear
convergence can be described as follows:
\beq\label{dfp-st-new}
\ba{rcl}
K_0^{\DFP} &\Def& \frac{4 n L}{\mu} \ln \frac{L}{\mu},
\ea
\eeq
and, for all $k \geq K_0^{\DFP}$, the new estimate
\eqref{dfp-new} takes the following form:
$$
\ba{rcl}
\lambda_k &\leq& \left( \frac{4 n L}{\mu k} \ln \frac{L}
{\mu} \right)^{k/2} \lambda_0 \quad (\;\refLE{dfp-st-new}\;
\lambda_0).
\ea
$$
Thus, compared to the old result, we have improved the
factor $\frac{L^2}{\mu^2}$ up to~$\frac{L} {\mu} \ln
\frac{L}{\mu}$. Interestingly enough, the ratio
between the old starting moments
\eqref{dfp-st-prev}, \eqref{bfgs-st-prev} of the superlinear
convergence of DFP and BFGS and the new ones
\eqref{dfp-st-new}, \eqref{bfgs-st-new} have remained the
same, $\frac{L}{\mu}$, although the both estimates have
been improved.

It is also interesting whether the results, obtained in this
paper, can be applied to \emph{limited-memory} quasi-Newton
methods such as L-BFGS \cite{LiuNocedal1989}. Unfortunately,
it seems like the answer is negative. The main problem is
that we cannot say anything interesting about just a
\emph{few} iterations of BFGS. Indeed, according to our main
result, after $k$ iterations of BFGS, the initial
residual is contracted by the factor of the form $[
\exp(\frac{n}{k} \ln \frac{L}{\mu}) - 1]^k$. For all
values $k \leq n \ln \frac{L}{\mu}$, this contraction factor
is in fact bigger than 1, so the result becomes useless.

\section{Conclusions}

We have presented a new theoretical analysis of local
superlinear convergence of classical quasi-Newton methods
from the convex Broyden class. Our analysis has been based
on the potential function involving the logarithm of
determinant of Hessian approximation and the trace of
inverse Hessian approximation. Compared to the previous
works, we have obtained new convergence rate estimates,
which have much better dependency on the condition number of
the problem.

Note that all our results are \emph{local}, i.e. they
are valid under the assumption that the starting point is
sufficiently close to a minimizer. In particular, there is
no contradiction between our results and the fact that
the DFP method is not known to be globally convergent with
inexact line search (see, e.g., \cite{ByrdNocedalYuan1987}).

Let us mention several open questions. First,
looking at the starting moment of superlinear convergence of
the BFGS method, in addition to the dimension of the
problem, we see the presence of the logarithm of its
condition number. Although typically such logarithmic
factors are considered small, it is still interesting to
understand whether this factor can be completely removed.

Second, all the superlinear convergence rates, which we have
obtained for the convex Broyden class in this paper, are
expressed in terms of the parameter~$\tau$, which controls
the weight of the DFP component in the updating formula for
the \emph{inverse} operator. At the same time,
in~\cite{RodomanovNesterov2020b}, the corresponding
estimates were presented in terms of the parameter $\phi$,
which controls the weight of the DFP component in the
updating formula for the \emph{primal} operator. Of course,
for the extreme members of the convex Broyden class, DFP and
BFGS, $\phi$ and $\tau$ coincide. However,
in general, they could be quite different. We do not know if
it is possible to express the results of this paper in terms
of $\phi$ instead of $\tau$.

Finally, in all the methods, which we considered, the
initial Hessian approximation $G_0$ was $L B$, where $L$ is
the Lipschitz constant of the gradient, measured relative to
the operator $B$. We always assume that this constant is
known. Of course, it is interesting to develop some
\emph{adaptive} algorithms, which could start from any
initial guess $L_0$ for the constant~$L$, and then somehow
dynamically adjust the Hessian approximations in iterations,
yet retaining all the original efficiency estimates.

\appendix

\section*{Appendix}
\setcounter{section}{1}

\BL
Let $A, G : \E \to \E^*$ be self-adjoint positive definite
linear operators, let $u \in \E$ be non-zero, and let
$\tau \in \R$ be such that $G_+ \Def \Broyd_{\tau}(A, G, u)$
is well-defined. Then,
\beq\label{inv-broyd}
\ba{rcl}
G_+^{-1} &=& \tau \left[ G^{-1} - \frac{G^{-1}
A u u^* A G^{-1}}{\la A G^{-1} A u, u \ra} + \frac{u
u^*}{\la A u, u \ra} \right] \\
&&\hspace{-1em} + \ (1 - \tau) \left[ G^{-1} - \frac{G^{-1}
A
u u^* + u
u^* A G^{-1}}{\la A u, u \ra} + \left( \frac{\la A G^
{-1} A u, u \ra}{\la A u, u \ra} + 1 \right) \frac{u u^*}
{\la A u, u \ra} \right],
\ea
\eeq
and
\beq\label{det-broyd}
\ba{rcl}
\Det(G_+^{-1}, G) &=& \tau \frac{\la A u, u \ra}{\la A
G^{-1} A u, u \ra} + (1 - \tau) \frac{\la G u, u \ra}{\la A
u, u \ra}.
\ea
\eeq
\EL

\begin{proof}
Denote $\phi \Def \phi_{\tau}(A, G, u)$. According
to Lemma~6.2 in \cite{RodomanovNesterov2020b}, we have
$$
\ba{rcl}
\Det(G^{-1}, G_+) &=& \phi \frac{\la A G^{-1} A u, u
\ra}{\la A u, u \ra} + (1 - \phi) \frac{\la A u, u
\ra}{\la G u, u \ra} \;\refEQ{def-phi}\; \left[ \tau 
\frac{\la A u, u \ra}{\la A G^{-1} A u, u \ra} + (1 - \tau)
\frac{\la G u, u \ra}{\la A u, u \ra} \right]^{-1}.
\ea
$$
This proves \eqref{det-broyd} since $\Det(G_+^{-1}, G) =
\frac{1}{\Det(G^{-1}, G_+)}$. Let us prove 
\eqref{inv-broyd}. Denote
\beq\label{def-G0-s}
\ba{rclrcl}
G_0 &\Def& G - \frac{G u u^* G}{\la G u, u \ra} + \frac{A u
u^* A}{\la A u, u \ra},
\qquad
s &\Def& \frac{A u}{\la A u, u \ra} - \frac{G u}{\la G u, u
\ra}.
\ea
\eeq
Note that
\beq\label{Gp-G0}
\ba{rcl}
G_+ &\refEQ{def-broyd}& G_0 + \phi \left[ \frac{\la G u, u
\ra A u u^* A}{\la A u, u \ra^2} + \frac{G u u^* G}{\la G u,
u \ra} - \frac{\la A u u^* G + G u u^* A}{\la A u, u \ra}
\right] \;=\; G_0 + \phi \la G u, u \ra s s^*.
\ea
\eeq
Let $I_{\E}$ and $I_{\E^*}$ be the identity operators in
$\E$ and $\E^*$. Since $G_0 u = A u$, we have
$$
\ba{rcl}
&& \left[ \left( I_{\E} - \frac{u u^* A}{\la A u, u \ra}
\right) G^{-1} \left(I_{\E^*} - \frac{A u u^*}{\la A u, u
\ra} \right) + \frac{u u^*}{\la A u, u \ra} \right] G_0 \\
&=& \left( I_{\E} - \frac{u u^* A}{\la A u, u \ra}
\right) G^{-1} \left( G_0 - \frac{A u u^* A}{\la A u, u
\ra} \right) + \frac{u u^* A}{\la A u, u \ra} \\
&\refEQ{def-G0-s}& \left( I_{\E} - \frac{u u^* A}{\la A u, u
\ra} \right) G^{-1} \left( G - \frac{G u u^* G} {\la G u, u
\ra} \right) + \frac{u u^* A}{\la A u, u \ra} \;=\; I_{\E}.
\ea
$$
Hence, we can conclude that
$$
\ba{rcl}
G_0^{-1} &=& \left( I_{\E} - \frac{u u^* A}{\la A u, u \ra}
\right) G^{-1} \left(I_{\E^*} - \frac{A u u^*}{\la A u, u
\ra} \right) + \frac{u u^*}{\la A u, u \ra} \\
&=& G^{-1} - \frac{G^{-1} A u u^* + u u^* A G^{-1}}{\la A
u, u \ra} + \left( \frac{\la A G^{-1} A u, u \ra}{\la A u, u
\ra} + 1 \right) \frac{u u^*}{\la A u, u \ra}.
\ea
$$
Thus, we see that the right-hand side of \eqref{inv-broyd}
equals
\beq\label{def-Hp}
\ba{rcl}
H_+ &\Def& G_0^{-1} - \tau \left[ \frac{\la A G^{-1} A u, u
\ra
u u^*}{\la A u, u \ra^2} + \frac{G^{-1} A u u^* A G^{-1}}
{\la A G^{-1} A u, u \ra} - \frac{G^{-1} A u u^* + u
u^* A G^{-1}} {\la A u, u \ra} \right] \\ &=& G_0^{-1} -
\tau \la A G^{-1} A u, u \ra w w^*,
\ea
\eeq
where
\beq\label{def-w}
\ba{rcl}
w &\Def& \frac{G^{-1} A u}{\la A G^{-1} A u, u \ra} -
\frac{u}{\la A u, u \ra}.
\ea
\eeq

It remains to verify that $H_+ G_+ = I_{\E}$. Clearly,
\beq\label{w-s}
\ba{rcl}
\la A G^{-1} A u, u \ra G_0 w &\refEQ{def-w}& G_0 G^{-1} A u
- \frac{\la A G^{-1} A u, u \ra G_0 u}{\la A u, u \ra} \\
&\refEQ{def-G0-s}& A u - \frac{\la A u, u \ra G u}{\la G
u, u \ra} \;\refEQ {def-G0-s}\; \la A u, u \ra s.
\ea
\eeq
Hence,
\beq\label{inv-aux}
\ba{rcl}
&&\la A G^{-1} A u, u \ra \la G_0 w, w \ra \;\refEQ
{w-s}\;
\la A u, u \ra \la s, w \ra \;\refEQ{def-w}\; \frac{\la A
u, u \ra \la s, G^{-1} A u \ra}{\la A G^{-1} A u, u \ra} -
\la s, u \ra \\
&& \qquad \;\refEQ{def-G0-s}\; \frac{\la A u, u \ra}{\la A
G^{-1} A u, u \ra} \left( \frac{\la A G^{-1} A u, u \ra}{\la
A u, u \ra} - \frac{\la A u, u \ra}{\la G u, u \ra} \right)
\;=\; 1 - \frac{\la A u, u \ra^2}{\la A G^{-1} A u, u \ra
\la G u, u \ra}.
\ea
\eeq
Consequently,
\beq\label{inv-aux2}
\ba{rcl}
\frac{\la G u, u \ra}{\la A u, u \ra} H_+ G_0 w w^* G_0
&\refEQ{def-Hp}& \frac{\la G u, u \ra}{\la A u, u \ra} (G_0^
{-1} - \tau \la A
G^{-1} A u, u \ra w w^*) G_0 w w^* G_0 \\
&=& \frac{\la G u, u \ra}{\la A u, u \ra} (1 -
\tau \la A G^{-1} A u, u \ra \la G_0 w, w \ra ) w w^* G_0 \\
&\refEQ{inv-aux}& \frac{\la G u, u \ra}{\la A u, u \ra}
\left( 1 - \tau + \tau \frac{\la A u, u \ra^2}{\la A G^{-1}
A u, u \ra \la G u, u \ra} \right) w w^* G_0 \\
&=& \left[ \tau \frac{\la A u, u \ra}{\la A G^{-1} A u, u
\ra} + (1 - \tau) \frac{\la G u, u \ra}{\la A u, u \ra}
\right] w w^* G_0.
\ea
\eeq
Thus,
\[
\ba[b]{rcl}
H_+ G_+ &\refEQ{Gp-G0}& H_+ ( G_0 + \phi \la G u, u \ra s
s^* ) \;\refEQ{w-s}\; H_+ \left( G_0 + \phi \frac{\la A
G^{-1} A u, u \ra^2}{\la A u, u \ra} \frac{\la G u, u \ra}
{\la A u, u \ra} G_0 w w^* G_0 \right) \\ &\refEQ{inv-aux2}&
H_+ G_0 + \phi \frac{\la A G^{-1} A u, u \ra^2}{\la A u, u
\ra} \left[ \tau \frac{\la A u, u \ra}{\la A G^{-1} A u, u
\ra} + (1 - \tau) \frac{\la G u, u \ra}{\la A u, u \ra}
\right] \\
&\refEQ{def-phi}& H_+ G_0 + \tau \la A G^{-1} A u, u \ra w
w^* G_0 \;\refEQ{def-Hp}\; I_{\E}.
\ea\qedEA
\]
\end{proof}

\begin{acknowledgements}
The presented results were supported by ERC Advanced Grant
788368. The authors are thankful to the anonymous reviewers
for their valuable time and comments.
\end{acknowledgements}


\begin{thebibliography}{10}

\bibitem{Davidon1959}
Davidon, W.: Variable metric method for minimization. Argonne
   National Laboratory Research and Development Report 5990
   (1959)

\bibitem{FletcherPowell1963}
Fletcher, R., Powell, M.: A rapidly convergent descent
   method for minimization. Computer Journal.
   \textbf{6}(2), 163--168 (1963)

\bibitem{Broyden1970p1}
Broyden, C.: The convergence of a class of double-rank
   minimization algorithms: 1.~General considerations.
   IMA Journal of Applied Mathematics. \textbf{6}(1),
   76--90 (1970)

\bibitem{Broyden1970p2}
Broyden, C.: The convergence of a class of double-rank
   minimization algorithms: 2.~The new algorithm. IMA
   Journal of Applied Mathematics. \textbf{6}(3), 222--231
   (1970)

\bibitem{Fletcher1970}
Fletcher, R.: A new approach to variable metric algorithms.
   Computer Journal. \textbf{13}(3), 317--322 (1970)

\bibitem{Goldfarb1970}
Goldfarb, D.: A family of variable-metric methods derived by
   variational means. Mathematics of Computation.
   \textbf{24}(109), 23--26 (1970)

\bibitem{Shanno1970}
Shanno, D.: Conditioning of quasi-Newton methods for function
   minimization. Mathematics of Computation.
   \textbf{24}(111), 647--656 (1970)

\bibitem{Broyden1967}
Broyden, C.: Quasi-Newton methods and their application to
   function minimization. Mathematics of Computation.
   \textbf{21}(99), 368--381 (1967)

\bibitem{DennisMore1977}
Dennis, J., Mor{\'e}, J.: Quasi-Newton methods, motivation
   and theory. SIAM Review. \textbf{19}(1), 46--89
   (1977)

\bibitem{NocedalWright2006}
Nocedal, J., Wright, S.: Numerical optimization.
   Springer Science \& Business Media, New York, NY,
   USA (2006)

\bibitem{LewisOverton2013}
Lewis, A., Overton, M.: Nonsmooth optimization via
   quasi-Newton methods. Mathematical Programming.
   \textbf{141}(1-2), 135--163 (2013)

\bibitem{Powell1971}
Powell, M.: On the convergence of the variable metric
   algorithm. IMA Journal of Applied Mathematics.
   \textbf{7}(1), 21--36 (1971)

\bibitem{Dixon1972p1}
Dixon, L.: Quasi-Newton algorithms generate identical points.
   Mathematical Programming. \textbf{2}(1), 383--387
   (1972)

\bibitem{Dixon1972p2}
Dixon, L.: Quasi Newton techniques generate identical points
   II: The proofs of four new theorems. Mathematical
   Programming. \textbf{3}(1), 345--358 (1972)

\bibitem{BroydenDennisMore1973}
Broyden, C., Dennis, J., Mor{\'e}, J.: On the local and
   superlinear convergence of quasi-Newton methods.
   IMA Journal of Applied Mathematics.
   \textbf{12}(3), 223--245 (1973)

\bibitem{DennisMore1974}
Dennis, J., Mor{\'e}, J.: A characterization of superlinear
   convergence and its application to quasi-Newton methods.
   Mathematics of Computation. \textbf{28}(126),
   549--560 (1974)

\bibitem{Stachurski1981}
Stachurski, A.: Superlinear convergence of Broyden's bounded
   $\theta$-class of methods. Mathematical
   Programming. \textbf{20}(1), 196--212 (1981)

\bibitem{GriewankToint1982}
Griewank, A., Toint, P.: Local convergence analysis for
   partitioned quasi-Newton updates. Numerische
   Mathematik. \textbf{39}(3), 429--448 (1982)

\bibitem{EngelsMartinez1991}
Engels, J., Mart{\'i}nez, H.: Local and superlinear
   convergence for partially known quasi-Newton methods.
   SIAM Journal on Optimization. \textbf{1}(1), 42--56
   (1991)

\bibitem{ByrdLiuNocedal1992}
Byrd, R., Liu, D., Nocedal, J.: On the behavior of
   Broyden's class of quasi-Newton methods. SIAM
   Journal on Optimization. \textbf{2}(4), 533--557 (1992)

\bibitem{YabeYamaki1996}
Yabe, H., Yamaki, N.: Local and superlinear convergence of
   structured quasi-Newton methods for nonlinear
   optimization. Journal of the Operations Research
   Society of Japan. \textbf{39}(4), 541--557 (1996)

\bibitem{WeiYuYuanLian2004}
Wei, Z., Yu, G., Yuan, G., Lian, Z.: The superlinear
   convergence of a modified BFGS-type method for
   unconstrained optimization. Computational
   Optimization and Applications. \textbf{29}(3), 315--332
   (2004)

\bibitem{YabeOgasawaraYoshino2007}
Yabe, H., Ogasawara, H., Yoshino, M.: Local and superlinear
   convergence of quasi-Newton methods based on modified
   secant conditions. Journal of Computational and
   Applied Mathematics. \textbf{205}(1), 617--632 (2007)

\bibitem{MokhtariEisenRibeiro2018}
Mokhtari, A., Eisen, M., Ribeiro, A.: IQN: An incremental
   quasi-Newton method with local superlinear convergence
   rate. SIAM Journal on Optimization.
   \textbf{28}(2), 1670--1698 (2018)

\bibitem{GaoGoldfarb2019}
Gao, W., Goldfarb, D.: Quasi-Newton methods: superlinear
   convergence without line searches for self-concordant
   functions. Optimization Methods and Software.
   \textbf{34}(1), 194--217 (2019)

\bibitem{RodomanovNesterov2020a}
Rodomanov, A., Nesterov, Y.: Greedy quasi-Newton methods with
   explicit superlinear convergence. CORE Discussion
   Papers. \textbf{06} (2020)

\bibitem{RodomanovNesterov2020b}
Rodomanov, A., Nesterov, Y.: Rates of Superlinear Convergence
   for Classical Quasi-Newton Methods. CORE Discussion
   Papers. \textbf{11} (2020)

\bibitem{JinMokhtari2020}
Jin, Q., Mokhtari, A.: Non-asymptotic Superlinear Convergence
   of Standard Quasi-Newton Methods. arXiv preprint
   arXiv:2003.13607 (2020)

\bibitem{ByrdNocedal1989}
Byrd, R., Nocedal, J.: A tool for the analysis of
   quasi-Newton methods with application to unconstrained
   minimization. SIAM Journal on Numerical Analysis.
   \textbf{26}(3), 727--739 (1989)

\bibitem{LiuNocedal1989}
Liu, D., Nocedal, J.: On the limited memory BFGS
   method for large scale optimization. Mathematical
   Programming. \textbf{45}(1-3), 503--528 (1989).

\bibitem{ByrdNocedalYuan1987}
Byrd, R., Nocedal, J., Yuan, Y.: Global convergence of a
   class of quasi-Newton methods on convex problems. SIAM
   Journal on Numerical Analysis. \textbf{24}(5), 1171--1190
   (1987).

\end{thebibliography}
\end{document}